# Fast and Multiscale Formation of Isogeometric matrices of Microstructured Geometric Models

T. Hirschler · P. Antolin · A. Buffa

**Abstract** The matrix formation associated to high-order discretizations is known to be numerically demanding. Based on the existing procedure of interpolation and lookup, we design a multiscale assembly procedure to reduce the exorbitant assembly time in the context of isogeometric linear elasticity of complex microstructured geometries modeled via spline compositions. The developed isogeometric approach involves a polynomial approximation occurring at the macro-scale and the use of lookup tables with pre-computed integrals incorporating the micro-scale information. We provide theoretical insights and numerical examples to investigate the performance of the procedure. The strategy turns out to be of great interest not only to form finite element operators but also to compute other quantities in a fast manner as for instance sensitivity analyses commonly used in design optimization.

**Keywords** Multiscale Mechanics · Matrix assembly · Isogeometric analysis · Additive Manufacturing · Geometric modeling · Lattice structures.

## 1 Introduction

The present work follows the paradigm for the design of microstructured geometries using functional composition as originally introduced by Elber (29). This approach has shown great flexibility to generate high-fidelity geometric models with heterogeneities as for instance lattice structures, structures made of composite materials or even random porous structures. More importantly, it initiates a flexible design framework of such structures since the geometries are naturally parameterized. However, in order to elaborate a complete methodology for the computational design of heterogeneous structures (or in other words microstructured structures), there is the need for efficient analysis methods adapted to these geometric models. Thanks to the emergence of the concept of IsoGeometric Analysis (IGA) introduced by Hughes et al. (41), it is now well-establish that spline-based geometric models present excellent performance for numerical simulations too (8; 17; 35; 47). Consequently, Massarwi et al. (53) and Antolin et al. (3) have been able to simulate the mechanical behavior of microstructured geometries by directly employing the geometric models developed in Elber (29). Despite the shown viability of the paradigm in Antolin et al. (3), there are several challenges to handle in order to build a practical method for designing complex heterogeneous geometries. The two main ones regarding the analysis arise from the high-order discretizations along with the explicit description of the heterogeneities that come with these geometric models. The numerical cost is tremendous if standard methods are blindly applied to the considered models.

Indeed, the presence of heterogeneities can be a source of complex local phenomena. Fine discretizations are often required to track them accurately in the context of the standard finite element method, which automatically leads to high computational cost. A multitude of alternative numerical approaches have been dedicated to multiscale problems to reduce the required computing resources. Methods based on model enrichment by means of augmented approximations spaces and superposition of solutions have shown great performances, particularly when the micro-scale is (drastically) smaller than the length scale of the structure. Among many others, let us mention the Multiscale FEM (MsFEM) (16; 19; 39) where numerically computed basis functions encode the micro-scale heterogeneities, the Generalized and Extended

T. Hirschler - P. Antolin - A. Buffa
Institute of Mathematics, Chair of Numerical Modelling and Simulation
École Polytechnique Fédérale de Lausanne, Switzerland

A. Buffa
Istituto di Matematica Applicata e Tecnologie Informatiche "Enrico Magenes"
Consiglio Nazionale delle Ricerche, Pavia, Italy

E-mail: thibaut.hirschler@epfl.ch

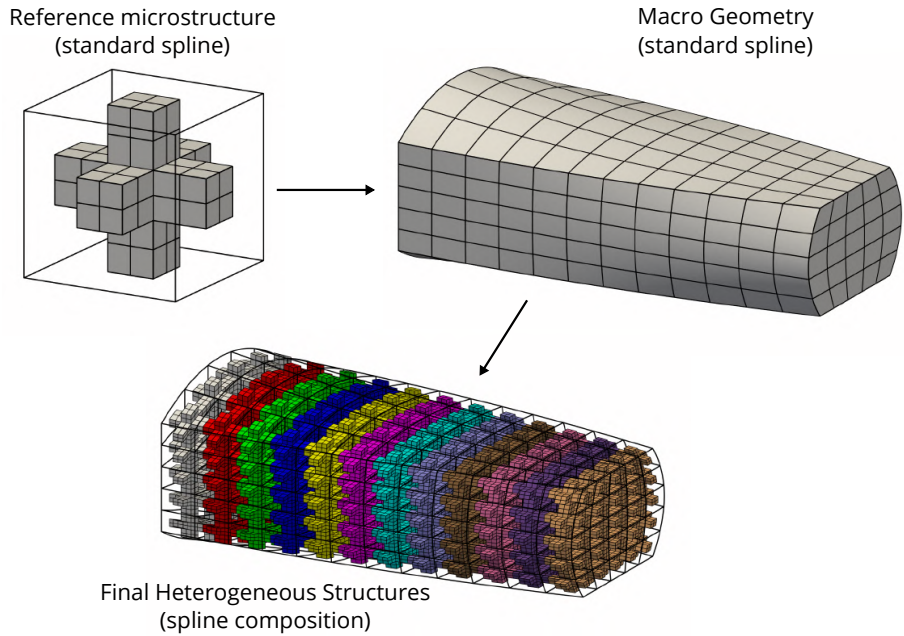

**Fig. 1** Geometric modeling of a complex heterogeneous structure using functional composition.

FMs (GFEM/XFEM) (26; 62; 67) based on partition of unity enrichment techniques, the multilevel FEM (FE$^2$) (33; 65) which employs higher-order homogenization scheme, or numerical homogenization (20; 56; 57) where the effective macroscopic description of the micro-heterogeneous material is identified by means of simulation on representative volume elements (RVE). Methods based on homogenization principle seems not applicable in our context. Here, the focus is on heterogeneous structures in the context of Additive Manufacturing and 3D printing (45). The precision of the printers constrains the achievable length scale range of the heterogeneities. We consider cases where the length scale of the heterogeneities is close to the length scale of the overall structures (factor 2 to 50). Let us also notice that the envisaged geometric models (29) provides the explicit representation of the microstructures. Therefore, we do not face the difficulty of mesh generation that motivates the use of some of the mentioned enhanced FEMs.

For this reason, we aim at solving the full scale problem and we concentrate on the use of fast assembly strategy in order to reduce the overall computational complexity. Instead of this first class of multiscale methods (superposition of micro/macro solutions), this work falls into a second class of numerical approaches for multiscale problems based on Domain Decomposition methods (13; 34; 43; 48). Instead of introducing the multiscale procedure during the modeling stage, here the idea is to derive specific multiscale solvers capable of solving the fine-scale finite element problem efficiently. Despite the resolution in itself, such a strategy is viable if and only if the formation of the full-field problem can be made in an acceptable amount of time. Our main source of inspiration are methods introduced to allow fast assembly of high-order IGA. In fact, this issue is known in the context of standard IGA (40) and successful procedures have been developed to reduce the assembly time in comparison with the standard finite element procedure based on element loop and Gauss quadrature. The main idea behind those fast formation procedures consists in taking advantage of the tensor-product nature of spline patches and the possible high inter-element continuity. More specifically, less costly quadrature rules have been developed as for instance reduced integration (31; 63) or specific rules (6; 18; 40; 42), and alternative formation procedures have been introduced as for instance sum factorization (2; 15) or the use of lookup tables (49; 50; 64). By combining several of these techniques, high order splines become tractable in isogeometric analysis (36; 58).

This leads us to develop a fast procedure to form the finite element operators associated to microstructured geometries modeled via spline compositions. To this purpose, we rely on the use of lookup tables as already investigated for standard IGA, and on the concept of unit cell inherited from multiscale strategies. The methodology is presented as follows. First, we present in Section 2 the employed geometric modeling approach based on spline composition. Section 3 starts by dealing with the simple case of the computation of the volume (Section 3.1), and then presents the developed assembly procedure (Section 3.2). Multiple numerical examples are given in Section 3.3 in order to study the performance of the strategy. Finally, Section 4 discusses the possible integration of the approach into a computational design framework for microstructured structures. Concluding remarks on our observations and findings are given in Section 5.

# 2 Geometric Modeling of Micro-structures

## 2.1 Multi-scale Modeling using Spline Composition

Functional composition is an attractive approach to generate geometric models of complex heterogeneous structures. The benefits of such an approach have been initially presented by Elber (29) and this modeling approach starts to gain visibility, see for instance Massarwi et al. (53); Antolin et al. (3). The main idea is depicted in Figure 1. It involves two principal ingredients:

– a macro-representation of the structure where the heterogeneities are not represented,
– a reference micro-structure which defines the pattern which will be tiled into the macro geometry.

Both geometric objects are modeled using standard CAD tools, as for instance splines. Without providing the theoretical details at the moment, let us mention that a discretization is resulting from the spline representation of these objects. Therefore, the macro-geometry can be viewed as a collection of non-overlapping macro-elements. The final heterogeneous structure, with an explicit representation of the micro-structure, is then obtained by embedding the reference micro-pattern into each of these macro-elements. This process is a composition: The overall geometric mapping that defines the geometry of the heterogeneous structure is prescribed by the composition of the mapping associated to reference micro-pattern and the mapping associated to the macro-representation of the structure.

It should be noticed that, as a consequence of the composition, the micro-pattern is not simply repeated into the macro-geometry but it is subjected to shape deformations too. In the case of the cross-tile depicted in Figure 1, it could be interesting to define the macro-mapping such that the arms of the microstructures are aligned with the directions associated to the principal stresses. This offers great flexibility, especially in a context of design optimization (3; 46). Futhermore, one can introduce more than one reference micro-structure. Different reference micro-patterns can be set to each macro-elements. All the examples tackled along this document highlight the great flexibility of the considered modeling approach. Finally, functional composition is not only interesting to model heterogeneous structures. For instance, it has been used to model cable reinforcements into membrane structures (7), or to model stiffened aerostructures (37).

## 2.2 Problem definition

### 2.2.1 Spline modeling

Splines are common tools in Computer Aided-Design as they are able to provide parameterization of general free-form objects. We employ, in this work, B-Splines and NURBS but the presented approach can be generalized to other spline technologies without any particular additional requirements. Mathematically, a spline is defined as a function of the form:

$$\mathscr{S} : \bar{\Omega}_S \to \Omega_S \,;\; \bar{x} \mapsto \mathscr{S}(\bar{x}), \qquad (1)$$

where the parameter $\bar{x}$ take values in the parameter space $\bar{\Omega}_S$ and is mapped to the codomain $\Omega_S$ which describes the domain delimited by the geometric object. In the case of a single B-Spline, the expression of $\mathscr{S}$ is given by:

$$\mathscr{S}(\bar{x}) = \sum_{i=1}^{n} R_i(\bar{x})\, \mathbf{s}_i \,, \qquad (2)$$

where the $n$ basis functions $R_i$ are piece-wise polynomial functions, which are weighted by the coefficients $\mathbf{s}_i$. We refer the reader to (22; 25; 32; 59) for the details on spline constructions. Here we just need to remind that splines are constructed on a non-uniform tensor product partition of the unit cube $\bar{\Omega}_S$. The lines and surfaces of such partition are called "parametric" knot lines while their images through $\mathscr{S}$ are called "physical" knot lines (or just knot lines) and constitute a partition of $\Omega_S$.

This simple spline model can be enriched in many ways: *e.g.*, geometries can be described as a collections of spline patches (multipatch geometries), or via NURBS in the case of conic sections. We refer the interested readers to the following non-exhaustive list of works which deal with key points regarding the modeling with splines: theoretical description and practical use of B-Spline and NURBS (22; 25; 32; 59), the issue of trimming procedures and boundary representation (1; 14; 21; 24; 52), and the generation of analysis-suitable geometric models for complex structures (4; 23; 51; 54).

### 2.2.2 Domain decomposition and Composition of mappings

Several mappings are involved in the geometric modeling approach presented in Section 2.1. The reader can refer to Figure 2 all along this section. We first consider a reference tile (or reference microstructure) which, in our case, is prescribed as a multi-patch spline model of the form:

$$\mathscr{T}^r(\theta) = \sum_{i=1}^{n_T} R_i^h(\theta)\, \mathbf{t}_i, \quad \theta \in \bar{\Omega}_T. \qquad (3)$$

The superscript $(\cdot)^h$ indicates that the basis functions are related to the fine scale (or micro-scale) of the problem.

Additionally, the model of the macro-geometry comes into play. In this work, we mainly consider single-patch models and we express these mappings as follows:

$$\mathscr{M}(\xi) = \sum_{i=1}^{n_M} R_i^H(\xi)\, \mathbf{m}_i, \quad \xi \in \bar{\Omega}_M. \qquad (4)$$

Analogously, the superscript $(\cdot)^H$ indicates that the basis functions are now related to the coarse scale (or macro-scale) of the problem. The knot lines of this macro-mapping decomposes

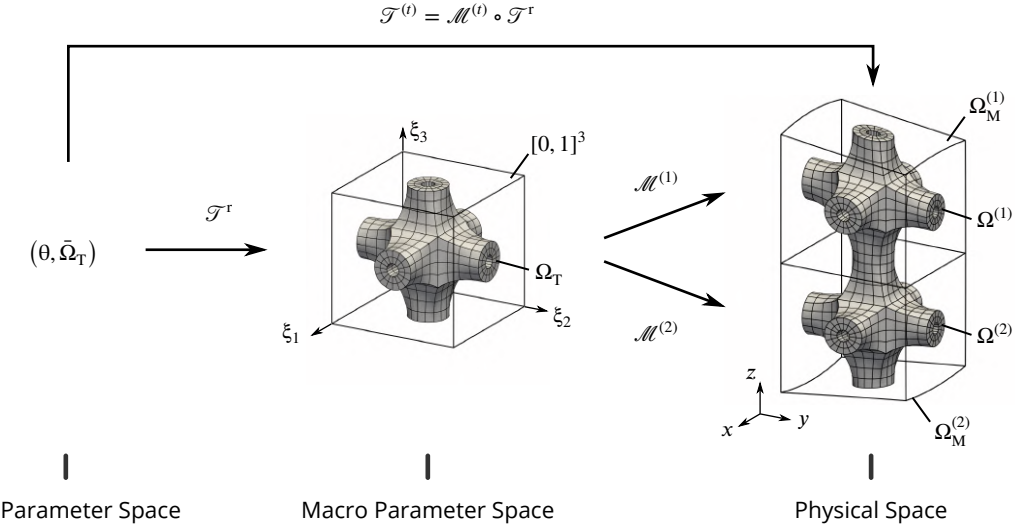

$$\mathcal{T}^{(t)} = \mathcal{M}^{(t)} \circ \mathcal{T}^{\mathrm{r}}$$

Fine Parameter Space · Macro Parameter Space · Physical Space

**Fig. 2** Definition of the spaces and the mappings involved in the multiscale geometric modeling of heterogeneous structures via functional composition.

the codomain $\Omega_{\mathrm{M}}$ into a non-overlapping domain decomposition (see Figure 3):

$$\Omega_{\mathrm{M}} = \bigcup_{e=1}^{m_{\mathrm{M}}} \Omega_{\mathrm{M}}^{(e)}, \tag{5}$$

where the subdomains are obtained by restricting the parameters to a subset $\bar{\Omega}_{\mathrm{M}}^{(e)}$ of the complete parameter space:

$$\Omega_{\mathrm{M}}^{(e)} = \left\{ \mathrm{x} \in \mathbb{R}^d \mid \mathrm{x} = \mathcal{M}(\xi), \; \forall \xi \in \bar{\Omega}_{\mathrm{M}}^{(e)} \right\}. \tag{6}$$

Even if it is not mandatory, let us rewrite these macro-subdomains through a Bézier extraction (11):

$$\Omega_{\mathrm{M}}^{(e)} = \left\{ \mathrm{x} \in \mathbb{R}^d \mid \mathrm{x} = \mathcal{M}^{(e)}(\xi), \; \forall \xi \in [0,1]^{\bar{d}} \right\}. \tag{7}$$

With this formulation, the parameter spaces associated to every sub-domains become identical and what differs from one macro-subdomain to the other is only the underlying geometric mapping. Each of them is equipped with its own Bézier mapping:

$$\mathcal{M}^{(e)}(\xi) = \sum_{i=1}^{n_{\mathrm{B}}} B_i^H(\xi) \, \mathbf{m}_i^{(e)}, \quad \xi \in [0,1]^{\bar{d}}. \tag{8}$$

where $B_i^H$ denotes the Bernstein basis polynomials. A Bézier mapping can be seen as a B-Spline with one single element. The Bézier extraction is a well-established procedure and takes the form of a linear application (*i.e.* the new coefficients $\mathbf{m}^{(e)}$ are expressed as a linear combination of the initial ones $\mathbf{m}$).

We now have in hand all the required ingredients to describe the fine-scale heterogeneous structure. The overall domain delimited by the full geometry is again given by a non-overlapping domain decomposition (see Figure 3):

$$\Omega = \bigcup_{t=1}^{m_{\mathrm{M}}} \Omega^{(t)}, \qquad \Omega^{(t)} = \left\{ \mathrm{x} \in \mathbb{R}^d \mid \mathrm{x} = \mathcal{T}^{(t)}(\theta), \; \forall \theta \in \bar{\Omega}_{\mathrm{T}} \right\}, \tag{9}$$

where the mappings $\mathcal{T}^{(t)}$ associated to the subdomains $\Omega^{(t)}$ are obtained through the composition of the reference tile (3) and the macro-elements (8):

$$\mathcal{T}^{(t)}(\theta) = \left( \mathcal{M}^{(t)} \circ \mathcal{T}^{\mathrm{r}} \right)(\theta) = \sum_{i=1}^{n_{\mathrm{B}}} B_i^H \left( \mathcal{T}^{\mathrm{r}}(\theta) \right) \mathbf{m}_i^{(t)}, \quad \theta \in \bar{\Omega}_{\mathrm{T}}. \tag{10}$$

The compositions are possible under the condition that the reference tile lies in the parameter space of the macro-elements, *i.e.* $\Omega_{\mathrm{T}} \subset [0,1]^{\bar{d}}$. Figure 2 summarizes the spaces and the mappings that are involved in the geometric modeling approach.

### 2.3 IsoGeometric Analysis
2.3.1 Test and trial functions

Following the isoparametric concept that is at the basis of Iso-Geometric Analysis (25; 41), the natural choices for trial and test functions are, for object $\mathcal{S}$, of the type:

$$u(\mathrm{x}) = \sum_{i=1}^{n} R_i \circ \mathcal{S}^{-1}(\mathrm{x}) \, \mathbf{u}_i, \tag{11}$$

where $R_i$ are the spline basis functions originated from the considered geometric model of the structure (see Equation (2)). The construction of the corresponding stiffness or mass matrices is cumbersome, due to the structure of these functions, especially in the case of the introduced geometric models which involve composition of mappings.

2.3.2 Challenges

The problem comes from the composition which may lead to a geometric mapping of *very* high degree. Indeed, it is possible to rewrite the composed map (10) as (multi-patch) tensor-product B-Splines (or NURBS) (69). However, if the reference microstructure $\mathscr{T}^r$ is a trivariate of degree $p_T$ per direction and the macro-geometry $\mathscr{M}$ is a trivariate of degree $p_M$ per direction, then the composed mapping is a trivariate of degree $p$ per direction, where $p$ is given by:

$$p = 3p_M p_T. \qquad (12)$$

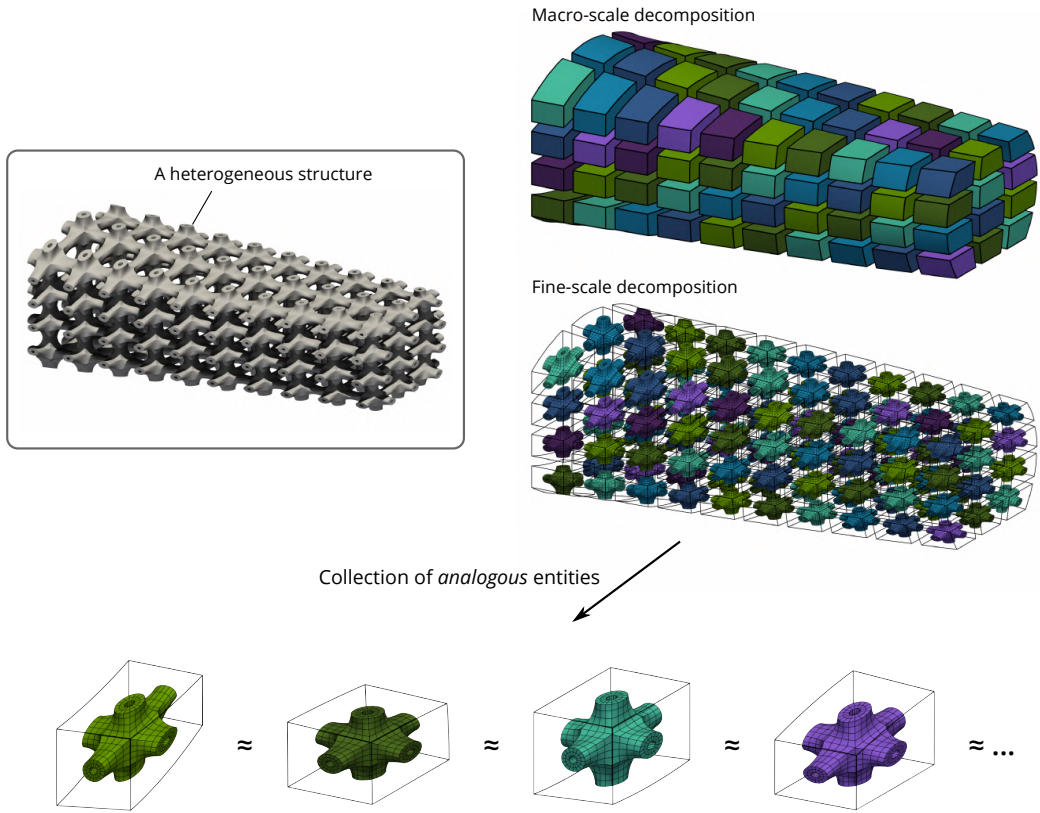

**Fig. 3** The considered geometric modeling approach naturally introduces a non-overlapping domain decomposition. Furthermore, the fine scale can be seen as a collection of microstructural entities that share many commonalities.

Thus, as an example, composing two cubic trivariates gives composed elements with degree 27. High degree is known to be a source of troubles during, for example, the numerical integration that is performed to build the system. Standard quadrature rules (and standard element-by-element assembly strategies) are rarely employed in practice when the degree is higher than 4: they become too expensive. In order to deal with large degree, more elaborate assembly strategies have been developed, see for example (2; 5; 18; 36; 40; 42; 49; 50; 58). The main idea of these approaches consists in exploiting the tensor-product nature of B-Splines and the inter-element regularity during the assembly. However, tracking degrees higher than 10 is still challenging. Especially here, the regularity of the composed map is limited by the regularity of the underlying reference micro-structural entity (i.e. the geometric map $\mathscr{T}^r$) which reduces the performances of the aforementioned fast-assembly strategies. Finally, even if it was efficient, one might question the pertinence of performing analyses with polyno-mial functions of degree 20 and more (e.g., due to bad conditioning and dense structure).

Consequently, we choose to limit the degree involved in the approximation subspace. A versatile way to build this space consists in recycling the basis functions $R_i^h$ involved in the geometric mapping of the reference micro-structure (remember Equation (3)). More specifically, the displacement field is approximated over each subdomain as follows:

$$u^{(t)}(\mathrm{x}) = \sum_{i=1}^{n_T} \left( R_i^h \circ \mathscr{T}^{(t)-1} \right)(\mathrm{x}) \, \mathbf{u}_i^{(t)}, \quad \mathrm{x} \in \Omega^{(t)}. \qquad (13)$$

Besides keeping the degree of the solution under control, this choice has also the advantage of discretizing the displacement field identically for each subdomain. This is not always true in the case of the isoparametric formulation (e.g. when the global mapping is a NURBS). With this choice, we also do not need to reformulate the composed maps with standard (high-degree) B-Splines. The reference micro-structures and the macro-elements are given separately in the data structure. Finally, let us mention that a similar choice regarding the solution spaces is done in Massarwi et al. (53).

However, the use of this non-isoparametric formulation is not enough, as such, to build an efficient simulation strategy for the targeted heterogeneous structures. Depending on the complexity of the reference microstructure and the number of macro-elements, the construction of the linear system of equations is still (very) demanding. As we said, in order to get an attractive assembly strategy, we employ ingredients from the fast assembly strategies that have been successfully developed in the standard isogeometric framework, but more importantly, we seek to exploit the inherent repetitiveness inside the heterogeneous structures as shown in Figure 3. Instead of performing blindly an assembly strategy, the key to success is to rely on the similarities between the subdomains in order to avoid computing the same quantities multiple times.

# 3 Fast multi-scale Assembly of Finite Element matrices

## 3.1 A preliminary illustrative example

Before tackling the case of the finite element operators, let us start with a more simple quantity in order to present the philosophy of the developed approach. To that purpose, let us try to compute the volume of a heterogeneous structure.

### 3.1.1 Multi-scale philosophy

The first step consists in invoking the non-overlapping domain decomposition of our geometric model (remember Section 2.2.2 and Figure 3). The overall volume is given by the sum of the local volumes:

$$\mathrm{V} = \sum_{t=1}^{m_{\mathrm{M}}} \mathrm{V}^{(t)}, \quad \text{with} \quad \mathrm{V}^{(t)} = \int_{\Omega^{(t)}} \mathrm{dV}.$$

The second step involves the two mappings defining the subdomain. More specifically, we pull back the mapping defining the macro-element and the mapping defining the reference microstructure. According to Figure 2, we go from the physical space back to the fine parameter space. These two successive changes of variables lead to:

$$\mathrm{V}^{(t)} = \int_{\bar{\Omega}_{\mathrm{T}}} |\det J_{\mathcal{T}^{\mathrm{r}}}| |\det J_{\mathcal{M}^{(t)}} \circ \mathcal{T}^{\mathrm{r}}| \, \mathrm{d}\theta, \quad (14)$$

where $\det J_{\mathcal{T}^{\mathrm{r}}}$ and $\det J_{\mathcal{M}^{(t)}}$ denote the Jacobian determinants of the mappings $\mathcal{T}^{\mathrm{r}}$ and $\mathcal{M}^{(t)}$, respectively. It is interesting to notice that in the expression of the local volume as given in Equation (14), there is only one quantity that differs from one subdomain to the other. This quantity is marked by the superscript identifying the current subdomain: it is the Jacobian determinant associated to the current macro-element, *i.e.* $\det J_{\mathcal{M}^{(t)}}$. Due to this quantity, one can hardly exploit the redundancies between the subdomains: One can, at best, store large precomputed quantities, for instance, the evaluation of $\det J_{\mathcal{T}^{\mathrm{r}}}$ at every integration point, but this does not represent major savings with the price of increasing the computer memory.

The role of the third step is to a certain extent to approximate the subdomain-dependent quantities (here the Jacobian determinant associated to the macro-elements) from the integrals. To do so, a projection step onto a spline space is introduced. Without giving the details at the moment, let us formulate the result of the projection:

$$\left(\Pi^H |\det J_{\mathcal{M}^{(t)}}|\right)(\xi) = \sum_{k=1}^{n_{\pi}} N_k^H(\xi) \, d_k^{H(t)}, \quad \xi \in [0,1]^{\bar{d}}. \quad (15)$$

The macro-Jacobian determinant is now prescribed as a linear combination of macro-basis functions $N_k^H$ weighted by some coefficient $d_k^{H(t)}$. Importantly, the same projection space is chosen for each subdomain. By substituting the projected quantities (15) into the expression (14) of the volume, and commuting the sum and the integral, we obtain:

$$\mathrm{V}^{(t)} \approx \sum_{k=1}^{n_{\pi}} \left( d_k^{H(t)} \int_{\bar{\Omega}_{\mathrm{T}}} |\det J_{\mathcal{T}^{\mathrm{r}}}| \left(N_k^H \circ \mathcal{T}^{\mathrm{r}}\right) \mathrm{d}\theta \right). \quad (16)$$

As a result, the volume of each subdomain can be computed through a dot product between two vectors:

$$\mathrm{V}^{(t)} \approx \mathbf{t}^h \cdot \mathbf{d}^{H(t)}, \quad (17)$$

where $\mathbf{d}^{H(t)} \in \mathbb{R}^{n_{\pi}}$ is the vector that collects the projection coefficients, and $\mathbf{t}^h \in \mathbb{R}^{n_{\pi}}$ is the vector that collects the integrals:

$$t_k^h = \int_{\bar{\Omega}_{\mathrm{T}}} |\det J_{\mathcal{T}^{\mathrm{r}}}| \left(N_k^H \circ \mathcal{T}^{\mathrm{r}}\right) \mathrm{d}\theta. \quad (18)$$

The vector $\mathbf{t}^h$ constitutes the so-called lookup table. Despite that the projection basis functions $N^H$ are defined at the macro-scale, we identify the lookup table with the fine-scale superscript $h$ because it contains all the information regarding the reference microstructure (the fine scale). Finally, with these quantities in hand, the overall volume is approximated by:

$$\mathrm{V} \approx \sum_{t=1}^{m_{\mathrm{M}}} \mathbf{t}^h \cdot \mathbf{d}^{H(t)} = \mathbf{t}^h \cdot \sum_{t=1}^{m_{\mathrm{M}}} \mathbf{d}^{H(t)}. \quad (19)$$

With this last expression, one can see that the two scales have been split. Interestingly, the exact same lookup-table is used for each subdomain. It needs to be computed only once, no matter the number of tiles there are in the heterogeneous structures. In practice, the lookup-table is computed in a preprocessing step and stored in a database. Figure 4 gives an overview of the approach.

### 3.1.2 Projection at the macro-level

Let us examine the projection step already introduced in Equation (15) in more detail. In this work, this key step is done through the application of a L2-projector occurring at the macro-level:

$$\Pi^H : L^2\left([0,1]^{\bar{d}}\right) \to Q_{p_{\pi}}\left([0,1]^{\bar{d}}\right), \quad (20)$$

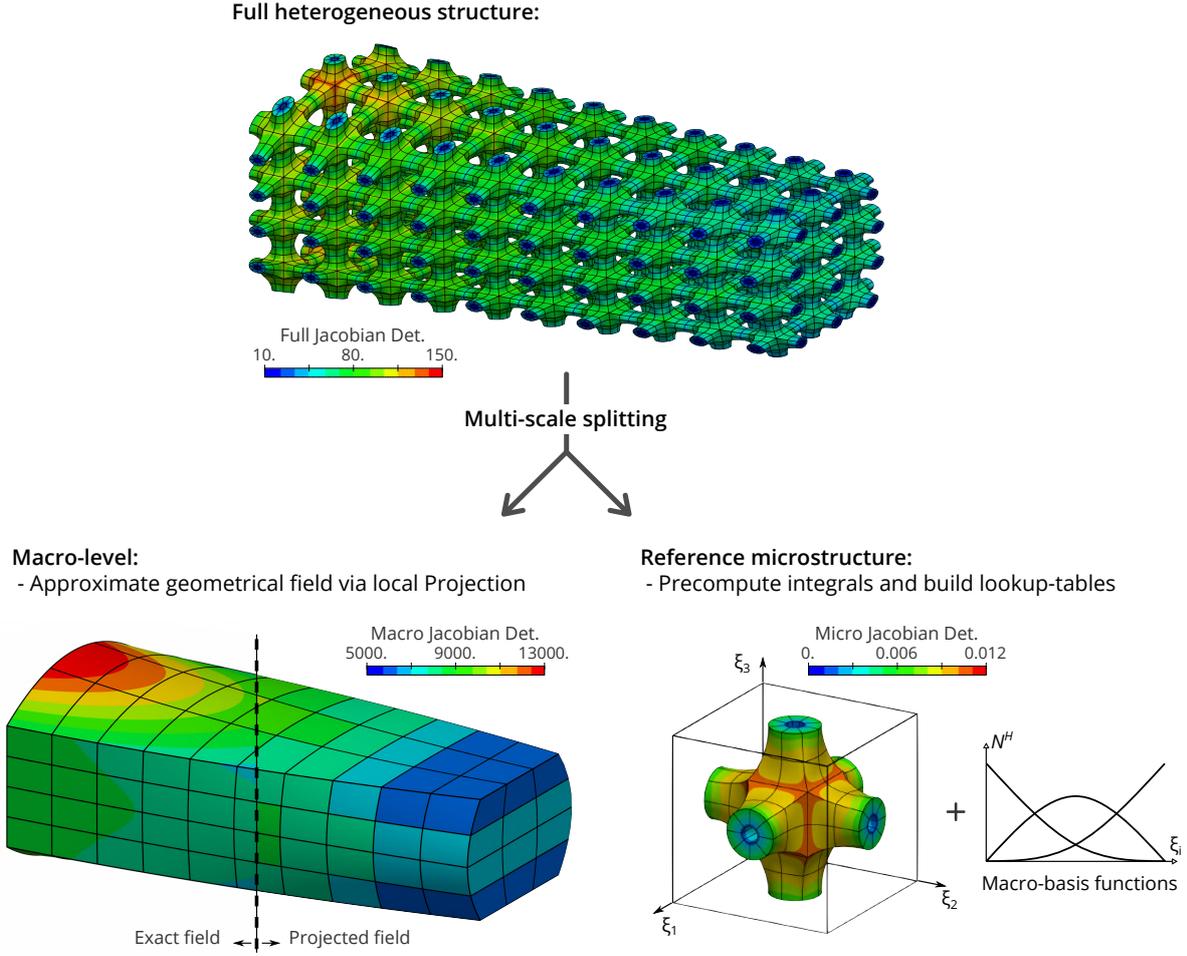

**Fig. 4** Main philosophy of the fast assembly approach for heterogeneous structures: complex fine-scale quantities are computed through the combination of projected macro-scale quantities and precomputed micro-scale quantities.

where the target subspace is a spline space of dimension $n_\pi = (p_\pi + 1)^{\bar{d}}$ defined as:

$$Q_{p_\pi} = \text{span}\left\{N_k^H\right\}_{k=1}^{n_\pi}, \qquad (21)$$

where the macro-basis function $N_k^H$ are taken as multivariate Bernstein basis polynomials of degree $p_\pi$. The projection may introduce some approximation, for instance, if the input function is not a polynomial or if the projection degree is lower than the degree of the input polynomial. For instance, if the global mapping is a B-Spline, then the Jacobian determinants associated to the macro-element are polynomials of degree:

$$\deg\left(\det J_{\mathcal{M}^{(t)}}\right) = \bar{d} p_M - 1. \qquad (22)$$

In order to get further numerical speed-ups, it can be interesting to balance between projection error and accuracy by taking the projection degree $p_\pi$ lower than the one of the Jacobian determinants. Furthermore, let us mention that in the case of NURBS, the Jacobian determinants might be rationals instead of polynomials. Consequently, we consider that, in general, some approximations are introduced during the projection step (either because the projection degree is lowered or because the input quantity is not a polynomial function).

In practice, the projection step mainly involves the resolution of a linear system of the form:

$$\mathbf{M}_\pi^H \mathbf{d}^{H(t)} = \mathbf{r}_\pi^{H(t)}, \qquad (23)$$

which enables to get the projection coefficients $\mathbf{d}^{H(t)}$. The components of the matrix $\mathbf{M}_\pi^H \in \mathbb{R}^{n_\pi \times n_\pi}$ and of the right-hand side $\mathbf{r}_\pi^{H(t)} \in \mathbb{R}^{n_\pi}$ defining the system to be solved are respectively given by:

$$\begin{aligned}\left[\mathbf{M}_\pi^H\right]_{kl} &= \int_{[0,1]^{\bar{d}}} N_k^H N_l^H \, d\xi, \\ \left(\mathbf{r}_\pi^{H(t)}\right)_k &= \int_{[0,1]^{\bar{d}}} N_k^H |\det J_{\mathcal{M}^{(t)}}| \, d\xi.\end{aligned} \qquad (24)$$

Selecting the same projection space for every subdomains leads to one unique projection matrix for each of the projection to be done. Thus, this matrix can be built and factorized (or even inverted) only once. In order to reduce the numerical effort even more, this matrix (or even its factorization/inverse) could be simply loaded from a database. In order to go even further, the resolution of the system (23) involved during the projection can be included into the lookup tables. Indeed, substituting Equation (23) into Equation (17) gives:

$$V^{(t)} \approx \mathbf{t}^h \cdot \mathbf{d}^{H(t)} = \mathbf{t}^h \cdot \mathbf{M}_\pi^{H-1} \mathbf{r}_\pi^{H(t)} = \mathbf{t}_\pi^h \cdot \mathbf{r}_\pi^{H(t)}. \quad (25)$$

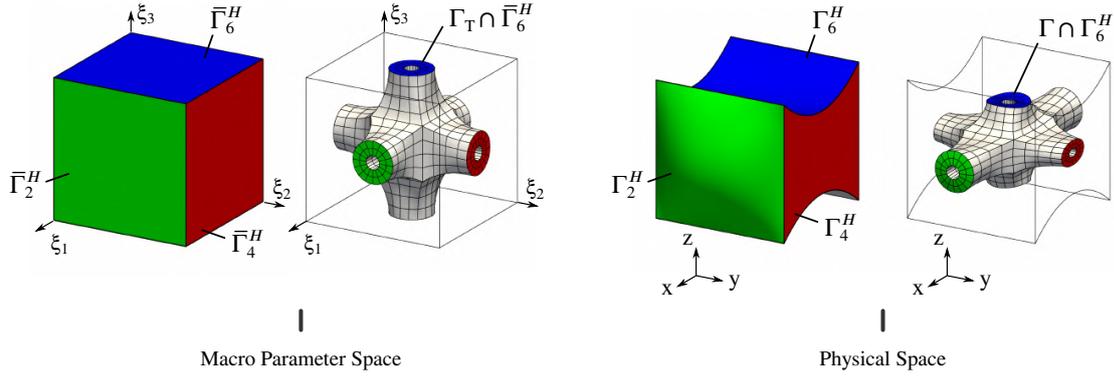

**Fig. 5** Definition of the boundaries where are prescribed external loads. Body forces are applied to the full solid body whereas the considered tractions/surface forces are applied to the colored surfaces.

The new lookup-table $\mathbf{t}_\pi^h$ now includes the inverse of the projection matrix:

$$\mathbf{t}_\pi^h = \mathbf{M}_\pi^{H-1} \mathbf{t}^h. \quad (26)$$

It is thus possible to somehow reduce the computational needs to the construction of the right-hand side $\mathbf{r}_\pi^{H(t)}$ only. However, even if it is a bit more computationally demanding, employing $\mathbf{t}^h$ instead of $\mathbf{t}_\pi^h$ has the advantage of computing the projection coefficients (by solving the system (23)). One can then directly evaluate the projected field using Equation (15), and therefore one can compute the projection error in a post-processing step.

### 3.1.3 Preliminary results

Let us consider the heterogeneous structure presented in Figure 4. Firstly, we compute the volume with a standard approach based on a numerical integration over the full fine-scale model. More specifically, we employ a Gauss-Legendre quadrature with 6 points per direction which leads, for this particular geometry, to compute the volume up to the machine precision. The computational time is about $t_1 \approx 6.5$ s. Secondly, we employ the developed fast strategy. In order to reach the machine precision too, quartic projection degree is required here (no error is introduced during the projection of the macro-Jacobian determinant with this choice). The computational time is more than 150 times lower, *i.e.* $t_2 \approx 0.04$ s (average value over multiple calculations) which represents a great computational speed-up. The construction of the lookup-table $\mathbf{t}^h$ is not taken into account in the computational time, and, in what follows, we always consider that the tables are built off-line, in a pre-processing step.

## 3.2 Application to Linear Elasticity

The strategy developed for the case of the volume computation in Section 3.1 can be applied analogously to the construction of common finite element operators. We detail in this section its application to linear elasticity.

### 3.2.1 Elastostatics

IGA relies on the same variational formulation than standard FEM. In the context of elastostatics, it consists in seeking the kinematically admissible displacement field $u \in \mathcal{U}$ such that:

$$W_{int}(u, v) + W_{ext}(v) = 0, \quad \forall v \in \mathcal{V}, \quad (27)$$

where $\mathcal{U}$ and $\mathcal{V}$ over the physical domain are add-hoc functional spaces that will contain the displacement solution and test functions, respectively. The internal virtual work $W_{int}$ and the external virtual work $W_{ext}$ in the equilibrium principle (27) are expressed as:

$$W_{int}(u, v) = -\int_\Omega \boldsymbol{\sigma}(u) : \boldsymbol{\varepsilon}(v) \, d\Omega, \quad (28)$$

$$W_{ext}(v) = \int_\Omega (v \cdot \mathbf{b}) \, d\Omega + \int_\Gamma (v \cdot \mathbf{t}) \, d\Gamma. \quad (29)$$

where $\boldsymbol{\varepsilon}$ and $\boldsymbol{\sigma}$ denote the linearized Green-Lagrange strain tensor and the linearized Cauchy stress tensor, respectively. The external loads are of two kinds: **b** denotes body forces (*e.g.* attraction or inertia forces), and **t** denotes tractions (or surface forces as a pressure force, for example) prescribed at the boundaries $\Gamma$ of the body. We distinguish two loading scenarios when dealing with tractions, as defined in Figure 5: the case where the load support is included in the macro-boundaries (colored faces in Figure 5), and the case where the load support lies inside the macro-elements (gray faces in Figure 5). This last case will not be considered in the scope of this work, as, in our practice, the loading scenarios and the material parameters are prescribed using the macro-model. Consequently, only tractions applied to the colored surfaces in Figure 5 can be described. The other scenario is automatically excluded.

The discretization of problem (27) that we consider is the Galerkin method applied on the functions constructed as in Equation (13). This choice generates a linear system:

$$\mathbf{K}\mathbf{u} = \mathbf{f}, \tag{30}$$

where $\mathbf{K}$ is the so-called stiffness matrix, $\mathbf{u}$ is the displacement vector that collects the Degrees Of Freedom (DOF), and $\mathbf{f}$ is the load vector. In what follows, we use the procedure proposed in Section 3.1 in order to perform fast assembly of the stiffness matrix and the load vector.

3.2.2 Polynomial approximation at the Macro-scale

*Pulling-back the macro-mappings* The non-overlapping domain decomposition enables to split the internal and external works into local contributions:

$$W_{\text{int}}(\boldsymbol{u}, \boldsymbol{v}) = \sum_{t=1}^{m_M} W_{\text{int}}^{(t)}(\boldsymbol{u}, \boldsymbol{v}), \quad W_{\text{ext}}(\boldsymbol{v}) = \sum_{t=1}^{m_M} W_{\text{ext}}^{(t)}(\boldsymbol{v}), \tag{31}$$

with:

$$W_{\text{int}}^{(t)}(\boldsymbol{u}, \boldsymbol{v}) = -\int_{\Omega^{(t)}} \boldsymbol{\sigma}(\boldsymbol{u}) : \boldsymbol{\varepsilon}(\boldsymbol{v}) \, \mathrm{d}\Omega, \tag{32}$$

$$W_{\text{ext}}^{(t)}(\boldsymbol{v}) = \int_{\Omega^{(t)}} (\boldsymbol{v} \cdot \mathbf{b}) \, \mathrm{d}\Omega + \int_{\Gamma^{(t)}} (\boldsymbol{v} \cdot \mathbf{t}) \, \mathrm{d}\Gamma. \tag{33}$$

As for the case of the volume, these local works are expressed at the level of the macro-parameter space after having pulled-back the macro-mappings $\mathscr{M}^{(t)}$. In order to identify the repeated quantities between the local works, one would need to reformulate them with respect to the curvilinear coordinates $\xi_i$ associated to the macro-elements (instead of the Cartesian coordinates associated to the physical space). Even if it is not mandatory, one can employ the curvilinear formalism commonly adopted for shell formulations (as presented in Bischoff et al. (9), for instance) to this purpose. We briefly detail this formalism in A.1 for the sake of clarity. Consequently, after pulling back the macro-mappings, we view the local works (32) and (33) as follows:

$$W_{\text{int}}^{(t)}(\boldsymbol{u}, \boldsymbol{v}) = -\int_{\Omega_T} \left( \boldsymbol{\sigma}^{(t)}(\boldsymbol{u}) : \boldsymbol{\varepsilon}^{(t)}(\boldsymbol{v}) \right) |\det J_{\mathscr{M}^{(t)}}| \, \mathrm{d}\xi, \tag{34}$$

$$W_{\text{ext}}^{(t)}(\boldsymbol{v}) = \int_{\Omega_T} \left( \boldsymbol{v} \cdot \bar{\mathbf{b}}^{(t)} \right) \mathrm{d}\xi + \int_{\Gamma_T} \left( \boldsymbol{v} \cdot \bar{\mathbf{t}}^{(t)} \right) \mathrm{d}\xi. \tag{35}$$

At this step, the involved quantities (stress tensor, strain tensor, and loads) are now marked with the sub-domain superscript because their expressions incorporate terms associated to the macro-mappings $\mathscr{M}^{(t)}$. Again, we refer to A.1 for the complete description of these quantities.

*Identification of the macro-fields to project* In contrast with the case of the volume, identifying the repeated quantities in the virtual works and thus the macro-quantities to project, is not straightforward, at least by looking at the expressions (34) and (35) only (especially for the internal work). Further manipulations of the integrands need to be performed, similarly to what is done for sum-factorization techniques, see for example (2) or (36). We detail in A.2 one way to rewrite the virtual works such that the quantities of interest (*i.e.* the sub-domain dependent quantities that will be projected) can be identified. Those developments leads us to view the local internal virtual works as:

$$W_{\text{int}}^{(t)}(\boldsymbol{u}, \boldsymbol{v}) = -\sum_{i=1}^{d} \sum_{j=1}^{d} \left( \int_{\Omega_T} \boldsymbol{u}_{,\xi_i}^{\top} \bar{\mathbf{A}}_{ij}^{H(t)} \boldsymbol{v}_{,\xi_j} \, \mathrm{d}\xi \right), \tag{36}$$

where the macro-fields denoted $\bar{\mathbf{A}}_{ij}^{H(t)}$ involve quantities associated to the macro-scale model only. These macro-fields are expressed by:

$$\begin{aligned}
\bar{\mathbf{A}}_{ij}^{H(t)} &: [0,1]^d \to \mathbb{R}^{d \times d}, \\
\xi &\mapsto \left[ \hat{\mathbf{G}}_i^{(t)} \check{\mathbf{C}}^{(t)} \hat{\mathbf{G}}_j^{(t)\top} \right] (\xi) \, |\det J_{\mathscr{M}^{(t)}}|(\xi),
\end{aligned} \tag{37}$$

with $d = 2, 3$ depending on the dimension of the problem. All the quantities involved in Equation (37) are described in detail in A.2. Importantly, the macro-fields do not only involve the geometrical information of the macro-model but they also include the material parameters. In this work, we consider isotropic elastic material, commonly defined through two parameters as for instance the Young's modulus $E$ and the Poisson's ratio $v$. We consider homogeneous material in the numerical examples in section 3.3 but heterogeneous material can also be modeled.

In contrast with the case of the volume where only one scalar-field is involved in the integral (recall Equation (14) and associated discussion), we end up here with a fourth-order tensor field. In the most general case, there are 45 distinct components for 3D problems (10 for 2D problems) due to following symmetries:

$$\bar{\mathbf{A}}_{ij}^{H(t)} = \bar{\mathbf{A}}_{ji}^{H(t)\top}. \tag{38}$$

For simple macro-geometries and homogeneous materials, the number of distinct components is then reduced. Figure 6 depicts the macro-fields involved in the internal work in the case of a free-form geometry. In general, these macro-fields are smooth rational polynomial functions defined over the macro-geometry. The polynomial degree of the denominator is the same than the one of the Jacobian determinant (recall Equation (22)). The numerator is also of high polynomial degree (50).

Regarding the external virtual work, the current multi-scale approach is applicable under the condition that the underlying quantities to project can be smoothly extended to the full macro-parameter space (as it was the case for the quantities involved in the internal work, *i.e.* the macro-fields are well defined over $[0,1]^d$ even if only their restrictions to $\Omega_T$ are used in the integrals). More specifically, the considered external works are expressed as:

$$W_{\text{ext}}^{(t)}(\boldsymbol{v}) = \int_{\Omega_T} \left( \boldsymbol{v} \cdot \bar{\mathbf{b}}^{H(t)} \right) d\xi + \int_{\Gamma_T} \left( \boldsymbol{v} \cdot \bar{\mathbf{t}}^{H(t)} \right) d\xi, \quad (39)$$

where $\bar{\mathbf{b}}^{H(t)}$ and $\bar{\mathbf{t}}^{H(t)}$ denote the extensions of the body and surface force terms, respectively. The body force term is define over the full macro-parameter space as follows:

$$\bar{\mathbf{b}}^{H(t)} : [0,1]^d \to \mathbb{R}^d; \xi \mapsto \bar{\mathbf{b}}^{H(t)}(\xi),$$
$$\text{such that} \quad \bar{\mathbf{b}}^{H(t)}|_{\Omega_T} = \left( \mathbf{b} \circ \mathcal{M}^{(t)} \right) |\det J_{\mathcal{M}^{(t)}}|. \quad (40)$$

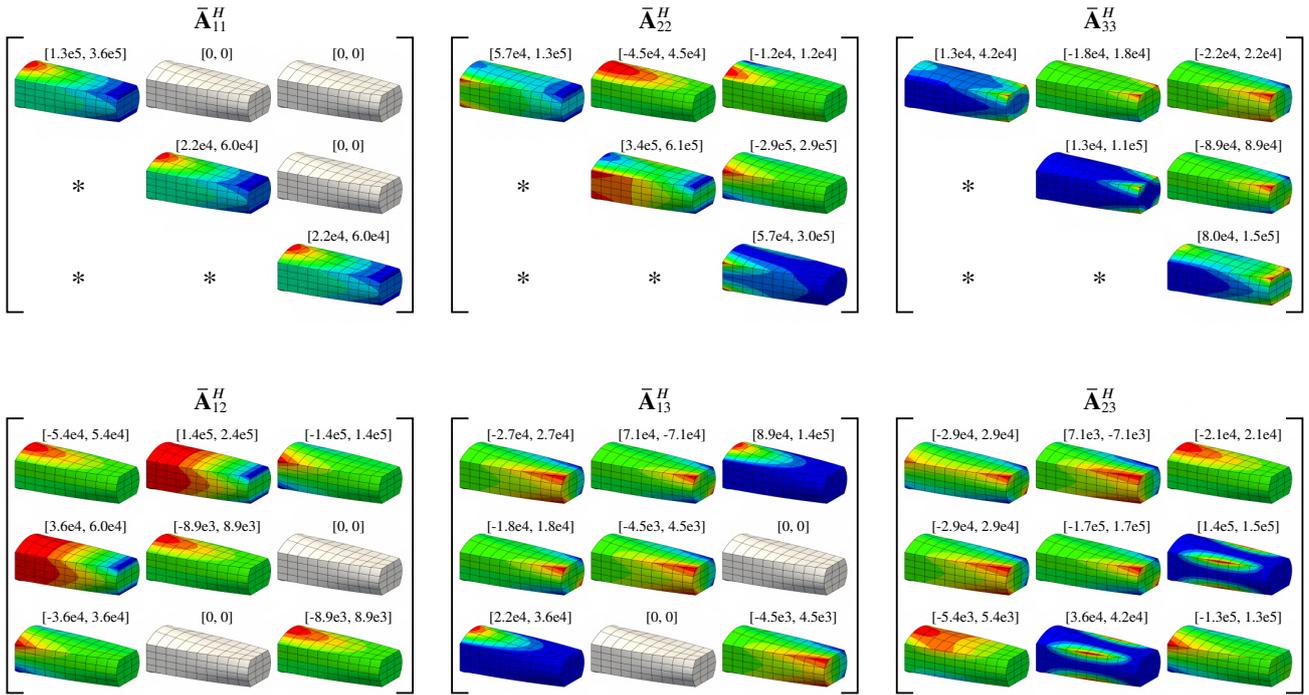

**Fig. 6** Macro field involved in the stiffness matrix

In the case of the surface forces, their extensions do not need to be defined over the complete macro-parameter space. Instead, a better choice seems to extend the associated terms over the macro-faces only:

$$\bar{\mathbf{t}}^{H(t)} : \bar{\Gamma}^H \to \mathbb{R}^d; \xi \mapsto \bar{\mathbf{t}}^{H(t)}(\xi),$$
$$\text{such that} \quad \bar{\mathbf{t}}^{H(t)}|_{\Gamma_T} = \left( \mathbf{t} \circ \mathcal{M}^{(t)} \right) |\det I_{\mathcal{M}^{(t)}}|. \quad (41)$$

By setting the boundary conditions using the macro-scale model as discussed in Section 3.2.1, we directly have in hand the required extensions of the load terms.

*Performing the projection* In order to get an exact spline representation of the macro-fields as defined in Equation (37), one would need to use a high degree NURBS (or rational Bézier). In this case, the required projection space might be different from one macro-element to the other (due to the weights) which is incompatible with our strategy: It breaks the possibility to build and use an unique lookup-table with pre-computed integrals that is common to every sub-domain. The required (possibly very) high degree to get the exact representation implies an important computational cost. Instead, we seek a balance between computational cost and precision: The macro-fields are projected on a spline space defined by multivariate Berstein polynomials of degree $p_\pi$ as for the volume case, see Equation (21). The projector defined in Equation (20) is applied to get a polynomial approximation of the components of macro-field involved in the internal work:

$$\Pi^H \bar{A}_{ij}^{H(t)}(\xi) = \sum_{C=1}^{n_\pi} N_C^H(\xi) \langle \bar{A}_{ij}^{H(t)} \rangle_C, \quad \xi \in [0,1]^d. \quad (42)$$

The projection coefficients marked as $\langle \bar{A}_{ij}^{H(t)} \rangle_C$ are obtained by solving systems of the same form as (23).

A similar polynomial approximation is introduced in the external virtual work (39). In order to deal with the surface forces, we introduce additional L2-projectors defined over each macro-face (or edges for 2D problem):

$$\Pi_F^H : L^2\left( \bar{\Gamma}_F^H \right) \to Q_{p_\pi}\left( \bar{\Gamma}_F^H \right), \quad F = 1, \ldots, 2d. \quad (43)$$

It is then possible to introduce the following polynomial approximations:

$$\Pi^H \bar{\mathbf{b}}^{H(t)}(\xi) = \sum_{C=1}^{n_\pi} N_C^H(\xi) \langle \bar{\mathbf{b}}^{H(t)} \rangle_C, \quad \xi \in [0,1]^d, \quad (44)$$

$$\Pi_F^H \bar{\mathbf{t}}_F^{H(t)}(\xi) = \sum_{C \in I_F} N_C^H(\xi) \langle \bar{\mathbf{t}}_F^{H(t)} \rangle_C, \quad \xi \in \bar{\Gamma}_F^H. \quad (45)$$

As for the other macro-quantities, the projection coefficients are obtained through the resolution of linear systems similar to Equation (23). It is however of reduced size for the field associated to the surface forces: Only the active projection basis functions at the current face, *i.e.* contained in the set denoted $I_F$ in Equation (45), are involved in the system to be solved.

*Consistency error term in Strang's first lemma* The impact of introducing polynomial approximations into the variational formulation has been theoretical studied, for instance, in (49). The projection error due to the polynomial approximation generates a consistency error that can be analyzed with Strang's first lemma. This enables to appropriately choose the projection degree $p_\pi$ such that the difference between the numerical solution and the original variational formulation is under control.

More precisely, Strang's first lemma (66) provides the following upper bound of the overall approximation error:

$$|u - u^{\pi h}|_{H^1(\Omega)} \leq \beta \Bigg( \sup_{w^h \in \mathcal{V}^h} \frac{|W_{\text{ext}}(w^h) - W_{\text{ext}}^\pi(w^h)|}{|w^h|_{H^1(\Omega)}} + \inf_{v^h \in \mathcal{V}^h} \bigg( |u - v^h|_{H^1(\Omega)} + \sup_{w^h \in \mathcal{V}^h} \frac{|W_{\text{int}}(v^h, w^h) - W_{\text{int}}^\pi(v^h, w^h)|}{|w^h|_{H^1(\Omega)}} \bigg) \Bigg), \quad (46)$$

where $\beta \in \mathbb{R}_+$ is a constant independent of the discretization. The polynomial projections (42), (44), and (45) influence the consistency errors (*i.e.* the two suprema of the upper bound). As shown in (49, Corollary 12), the approximate bilinear form $W_{\text{int}}^\pi$ satisfies the inequality:

$$\sup_{w^h \in \mathcal{V}^h} \left( \frac{|W_{\text{int}}(v^h, w^h) - W_{\text{int}}^\pi(v^h, w^h)|}{|w^h|_{H^1(\Omega)}} \right) \leq \beta_2 |v^h|_{H^1(\Omega)} \varepsilon_\pi, \quad (47)$$

where $\varepsilon_\pi$ is the projection error given by:

$$\varepsilon_\pi = \max_{t=1, m_M} \left\{ \max_{\xi \in \Omega_T} \max_{ijkl=1,d} |\bar{A}_{ijkl}^{H(t)}(\xi) - \Pi^H \bar{A}_{ijkl}^{H(t)}(\xi)| \right\}. \quad (48)$$

The constant $\beta_2 \in \mathbb{R}_+$ takes into account the equivalence between the seminorms on the parameter domain and on the physical space (8; 49), and the inequality (47) is simply a Young's inequality.

From Equation (47), it can be inferred that the consistency error is kept under control during the projection. Indeed, we can evaluate the projection error in Equation (47) since we know the quantities we are projecting. This can be done by relying on a simple sampling-based approach where the differences between the exact macro-quantities and their projections are evaluated at each sample point. More importantly, this is done *a priori* to the assembly of the finite element operators. Consequently, a simple approach for choosing the projection degree $p_\pi$ consists in prescribing a given tolerance to the projection error (48) and an initial projection space. The projection degree is then increased recursively while the projection error gets below the given tolerance. A convenient choice for the initial guess is to set the projection degree equal to the spline degree of the macro-mappings. Regarding the prescribed tolerance, it might not be necessary to set it too low for two principal reasons. First, the discretization error might be predominant and reducing the consistency error below the discretization error does not improve the overall approximation error (see again Strang's first lemma (46)). Second, the computational cost for building the finite element operators is directly linked to the chosen projection space. An inappropriately too low tolerance on the projection error will lead to (possibly large) extra cost during the assembly. The numerical studies in Section 3.3 enable to further understand the link between the projection error, the overall solution error, and the assembly cost.

### 3.2.3 Forming lookup tables at the fine scale

Let us now elaborate a fast-assembly strategy of the stiffness matrix and of the load vector based on the strategy previously illustrated in Figure 4. At the current step, we know what are the macro-fields to project. However, we still need to define and build the so-called lookup tables as done for the case of the volume. These lookup tables are the same to every microstructural entity under the condition that the finite dimensional subspaces discretize the solution identically at the sub-domain level. The choice associated to Equation (13) fulfills this requirement.

The overall operators are obtained by assembling sub-operators associated to the sub-domains. By invoking the expression of the internal virtual work as given by Equation (36), the contribution to the local stiffness matrices associated to the test basis function $R_A^h$ and trial basis function $R_B^h$ reads as:

$$\mathbf{K}_{AB}^{(t)} = \sum_{i=1}^{d} \sum_{j=1}^{d} \int_{\Omega_T} R_{A',\xi_i}^h R_{B',\xi_j}^h \bar{A}_{ij}^{H(t)} \, d\xi . \quad (49)$$

By introducing the spline approximation of the macro-fields as provided in Equation (42) into this last expression, we get:

$$\mathbf{K}_{AB}^{\pi(t)} = \sum_{i=1}^{d} \sum_{j=1}^{d} \sum_{C=1}^{n_\pi} \langle \bar{A}_{ij}^{H(t)} \rangle_C \underbrace{\int_{\Omega_T} R_{A',\xi_i}^h R_{B',\xi_j}^h N_C^H \, d\xi}_{\text{Independent of } \mathcal{M}^{(t)}} . \quad (50)$$

The remaining integrals are independent of the macro-geometry and of the material parameters. The integrands only depend on the geometric model of the reference microstructure and the macro-projection space. Therefore, they can be precomputed and stored in so-called lookup tables:

$$\mathfrak{K}_r^h \in \mathbb{R}^{n_T \times n_T \times n_\pi \times d \times d}, \tag{51}$$

which takes the form of elementary stiffness matrices associated to the micro-geometry only, times some macro-basis polynomial functions associated to the projection space. The detailed expressions of the components are:

$$\mathfrak{K}_r^h(A, B, C, i, j) = \int_{\Omega_T} \Big\{ \ldots \\ \frac{d}{d\xi_i}\big(R_A^h \circ \mathcal{T}^{r-1}\big)(\xi) \frac{d}{d\xi_j}\big(R_B^h \circ \mathcal{T}^{r-1}\big)(\xi) N_C^H(\xi) \Big\} d\xi. \tag{52}$$

These integrals cannot be computed analytically, in general, due to the contribution of the geometric map describing the reference micro-structure. Instead, they are evaluated via an accurate enough numerical integration.

An identical procedure is followed in order to build the load vector. The contribution to the local load vector associated to the test basis function $R_A^h$ is approximated by:

$$\mathbf{f}_A^{\pi(t)} = \sum_{C=1}^{n_\pi} \langle \bar{\mathbf{b}}^{H(t)} \rangle_C \underbrace{\int_{\Omega_T} R_A^h N_C^H \, d\xi}_{\text{Independent of } \mathcal{M}^{(t)}} \\ + \sum_{F=1}^{2d} \sum_{C \in I_F} \langle \bar{\mathbf{t}}_F^{H(t)} \rangle_C \underbrace{\int_{\Gamma_T \cap \bar{\Gamma}_F^H} R_A^h N_C^H \, d\xi}_{\text{Independent of } \mathcal{M}^{(t)}}. \tag{53}$$

Again, we end with integrals that are independent of the macro-geometry and of the boundary conditions (Neumann type). These integrals are concatenated into a lookup table defined as:

$$\mathfrak{F}_r^h \in \mathbb{R}^{n_T \times n_\pi \times (1+2d)}, \tag{54}$$

where the detailed expressions of the components are:

$$\mathfrak{F}_r^h(A, C, 0) = \int_{\Omega_T} \big(R_A^h \circ \mathcal{T}^{r-1}\big)(\xi) N_C^H(\xi) \, d\xi, \\ \mathfrak{F}_r^h(A, C, F) = \int_{\Gamma_T \cap \bar{\Gamma}_F^H} \big(R_A^h \circ \mathcal{T}^{r-1}\big)(\xi) N_C^H(\xi) \, d\Gamma_T, \tag{55}$$

with $F \in [\![1, 2d]\!]$. Finally, as for the lookup table used for the construction of the stiffness matrix, the lookup table for the load vector is considered as precomputed and stored in a database.

*Implementation aspects* The lookup tables should not be built in practice as dense tensors since a lot of the components may be zero. More specifically, the non-zero components that need to be stored in the lookup table $\mathfrak{K}_r^h$ can be identified by the following set of basis function couples:

$$I_{\text{coo}} = \Big\{ (A, B) \in [\![1, n_T]\!]^2 \mid \\ A \leq B, \operatorname{supp} R_A^h \cap \operatorname{supp} R_B^h \neq \varnothing \Big\}. \tag{56}$$

The condition of selecting index $A$ lower than index $B$ is a way of taking advantage of the symmetry of the stiffness matrix to be built: Only the upper part is constructed (at least in a first step). Consequently, the lookup-tables are viewed as matrix objects:

$$\operatorname{mat}\big[\,\mathfrak{K}_r^h\,\big] \in \mathbb{R}^{n_{nz} \times n_\pi d^2}, \quad \text{with} \quad n_{nz} = \operatorname{card} I_{\text{coo}}. \tag{57}$$

Finally, all the non-zero components of the full stiffness matrix can be computed by performing a single matrix-matrix product:

$$\underbrace{\operatorname{nzdata}\big[\,\mathbf{K}\,\big]}_{(n_{nz} \times d^2 m_M)} = \underbrace{\operatorname{mat}\big[\,\mathfrak{K}_r^h\,\big]}_{(n_{nz} \times n_\pi d^2)} \cdot \underbrace{\operatorname{mat}\big[\,\langle \bar{\mathbf{A}}^H \rangle\,\big]}_{(n_\pi d^2 \times d^2 m_M)}, \tag{58}$$

where all the projection coefficients have been concatenated into a rectangular matrix. Once the non-zero values of the stiffness matrix in hand, one can employ any scientific library that handles sparse matrix data structures to generate the overall matrix object.

#### 3.2.4 The multiscale fast assembling strategy

Finally, let us summarize the main steps involved during the fast multi-scale assembly strategy:

1. perform the L2-projection of the macro-fields (Equations (42), (44), and (45)),
2. load the lookup tables associated to the reference microstructure (Equations (52) and (55)),
3. compute the non-zero values of the stiffness matrix and of the load vector (Equation (58)),
4. build the overall finite element operators using a sparse data structure.

The projections occur at the macro-scale only and therefore do not form the most expensive step. Step 3 is largely the most demanding step as it involves a matrix-matrix product of quite large matrices.

### 3.3 Numerical Examples

In order to analyze the performance of the presented multiscale assembly approach, we compute several quantities along the numerical examples. We compare the results obtained with this new fast assembly strategy against the very standard procedure based on element loops with full Gauss quadrature (36). More precisely, we define:

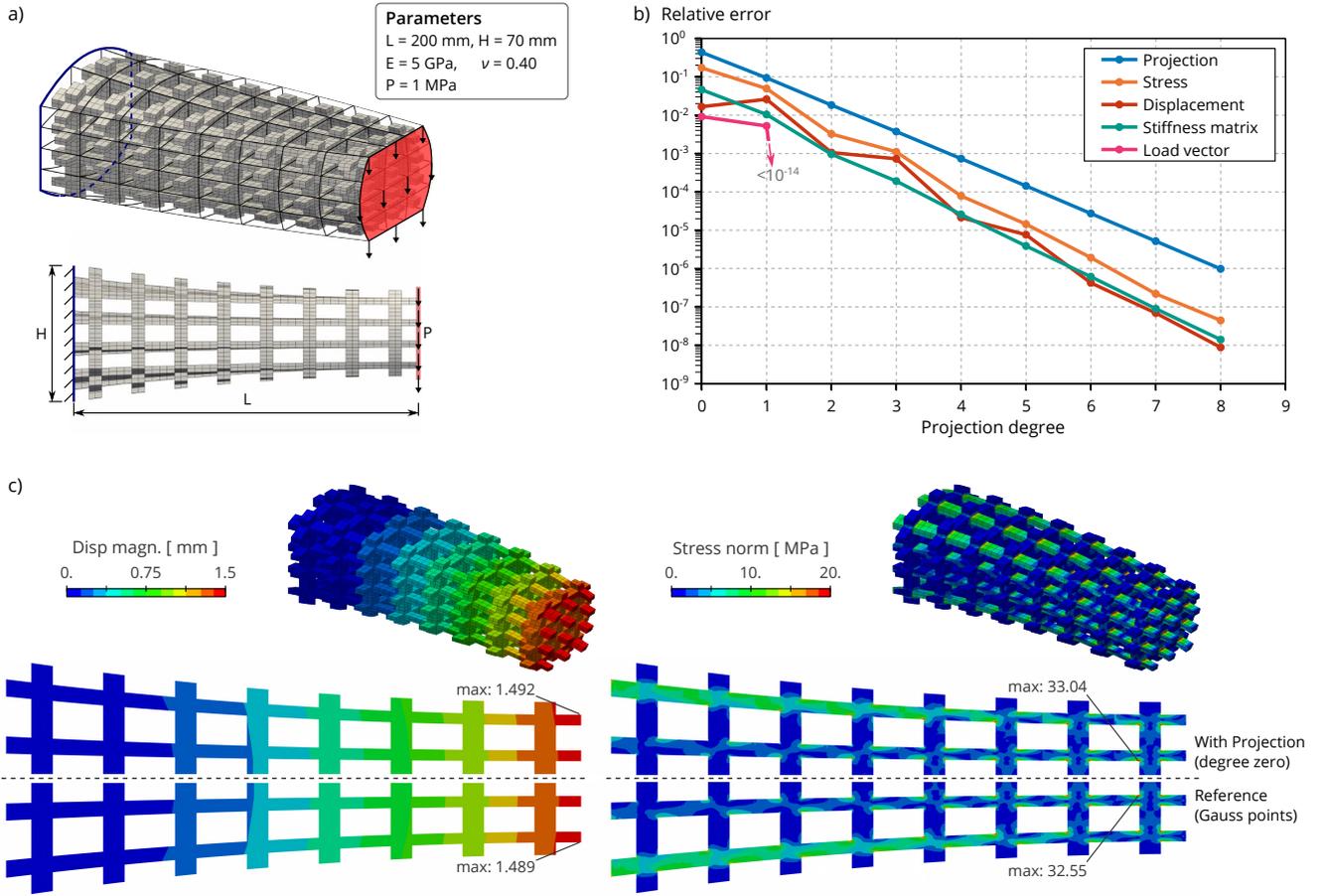

**Fig. 7** Influence of the projection space: We compare the numerical solution obtained with the operators involving different projection spaces with the numerical solution obtained with the standard operators via Gauss quadrature. a) Mechanical settings, b) Relative error quantities as defined in Equations (59)-(63), and c) Displacement and Stress fields.

– the relative projection error:

$$E_{\text{proj}} = \max_{t=1,m_M} \max_{\Omega_T} \|\bar{\mathbf{A}}^{H(t)} - \Pi^H \bar{\mathbf{A}}^{H(t)}\|_F / \|\bar{\mathbf{A}}^{H(t)}\|_F \quad (59)$$

– the relative solution errors (displacement and stress fields):

$$E_{\text{disp}} = |\boldsymbol{u}^\bullet - \boldsymbol{u}^\pi|_{H^1(\Omega)} / |\boldsymbol{u}^\bullet|_{H^1(\Omega)}, \quad (60)$$

$$E_{\text{stress}} = \|\boldsymbol{\sigma}^\bullet - \boldsymbol{\sigma}^\pi\|_{L^2(\Omega)} / \|\boldsymbol{\sigma}^\bullet\|_{L^2(\Omega)}. \quad (61)$$

– the stiffness matrix and load vector differences:

$$E_{\text{mat}} = \|\mathbf{K}^\bullet - \mathbf{K}^\pi\|_F / \|\mathbf{K}^\bullet\|_F, \quad (62)$$

$$E_{\text{vec}} = \|\mathbf{f}^\bullet - \mathbf{f}^\pi\|_F / \|\mathbf{f}^\bullet\|_F. \quad (63)$$

– the speed-up that evaluates the computational savings:

$$T_{\text{cost}} = \text{time}\left(\mathbf{K}^\bullet, \mathbf{f}^\bullet\right) / \text{time}\left(\mathbf{K}^\pi, \mathbf{f}^\pi\right). \quad (64)$$

*Note:* the provided assembly times in what follows measure only the time spent to compute the non-zero components and not the formation of the full matrix. Regarding the new multiscale procedure, we consider the lookup tables as pre-built. All the numerical experiments were computed on a 1.90GHz i7-8665U (Intel Core) processor with 32 GB RAM (serial computing). More importantly, the code that enables both standard and new assembly procedures has been developed by one single person into the same coding environment.

The superscripts $(\cdot)^\bullet$ and $(\cdot)^\pi$ indicate that the underlying quantity is obtained with the standard and the fast procedures, respectively. Considering the solutions, the results obtained with the standard procedure are taken as the reference solutions: We aim at obtaining results similar to those given by the standard procedure. Considering the computational time, we expect that the fast procedure reduces the time needed to assemble the linear systems (speed-up greater than one).

3.3.1 Influence of the projection error

Following the discussion associated to Equation (47), we know that there is a direct link between the projection error and the consistency error. When comparing the standard procedure based on Gauss quadrature and the new multiscale procedure, the projection error is directly linked to the solution errors in the same manner (the first term representing the discretization error vanishes when using $\boldsymbol{u} = \boldsymbol{u}^\bullet$ in Strang's lemma (46)). This theoretical behavior is observed numerically as shown in Figure 7. The reduction of the projection error implies the decrease of all the computed errors. It can also be seen that the solution obtained with the fast operators is very close to the one obtained with the standard operators. Even for a projection degree equal to zero ($p_\pi = 0$), the displacement field $\boldsymbol{u}^\pi$ and the stress field $\boldsymbol{\sigma}^\pi$ are quasi indistinguishable from the reference fields obtained via the standard procedure (with Gauss quadrature), see Figure 7 again. For higher projection degree, the differences between the solutions are not clearly visible, as reflected by the low relative errors.

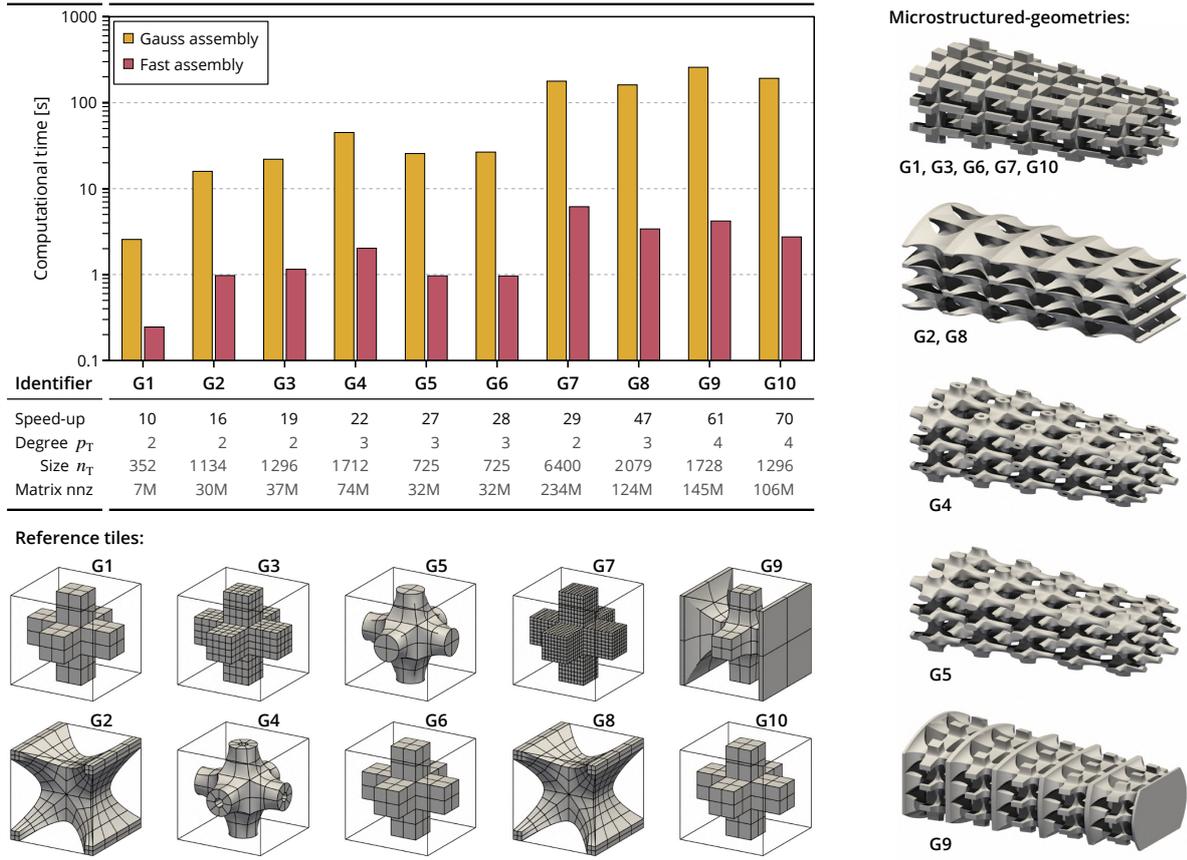

**Fig. 8** Influence of the reference tile: We study here the assembly times for ten different reference microstructures.

It is however not necessary to select a too high projection degree: Recall that the reference numerical solution (in our case $\boldsymbol{u}^\bullet$) might be far from the exact solution due to the discretization error. Finally, a practical strategy for selecting the projection degree $p_\pi$ consists in setting a maximal tolerance on the projection error. The selected projection degree is a solution of:

$$\text{Find} \quad p_\pi \in \mathbb{N} \quad \text{such that} \quad E_{\text{proj}}(p_\pi) \leq \text{tol}. \tag{65}$$

In what follows, we fix the tolerance on the projection error as tol $= 10^{-3}$. With such a choice, the required projection degree is $p_\pi = 4$ for the microstrucure tackled in Figure 7, which leads to numerical solution errors $E_{\text{disp}} = 2 \cdot 10^{-5}$ and $E_{\text{stress}} = 8 \cdot 10^{-5}$.

3.3.2 Influence of the reference microstructrure

Figure 8 summarizes the study of the influence of the reference microstructure on the computational time required to form the operators. Here, the macro-geometry is the same for each case and we vary the reference tile. Thus, the projection space is the same for each case (setting the tolerance on the projection error to $10^{-3}$ leads to $p_\pi = 4$ again). Ten different reference tiles with different spline degree $p_T$, mesh density $n_T$, or shape $\Omega_T$ are successively used. More precisely, we report in Figure 8 the computational time to build the operators with either the standard procedure based on Gauss quadrature and the presented fast procedure based on projection and lookup tables. We also give the resulting speed-up and some information regarding the underlying discretization (the spline degree $p_T$, the number of control points $n_T$, and the number of non-zero components in the full stiffness matrix).

By comparing the cases marked by the identifiers G1, G6, and G10, or by comparing the cases marked by the identifiers G2 and G8, one can observe the influence of the solution degree $p_T$. These reference tiles are identical in their shape and their number of element but are discretized with different spline degrees (quadratic, cubic, or quartic). It is not surprising that the speed-up increases with the solution degree. As mentioned, the presented approach originates from the assembly procedure of Mantzaflaris and Jüttler (49). This approach shows great benefits for high degree parametrization (49; 58). It is also known that the standard procedure based on element loop and Gauss quadrature is too expensive when it comes to high degree, especially for 3D problems (it is rarely employed for quartic spline degree and above). Therefore, the use of lookup tables with precomputed integrals gives access to high spline degree. Indeed, the assembly time for the case G10 (quartic solution space) is lower than 3 seconds with the fast approach, whereas it takes more than 3 minutes to assemble the operators with the standard procedure.

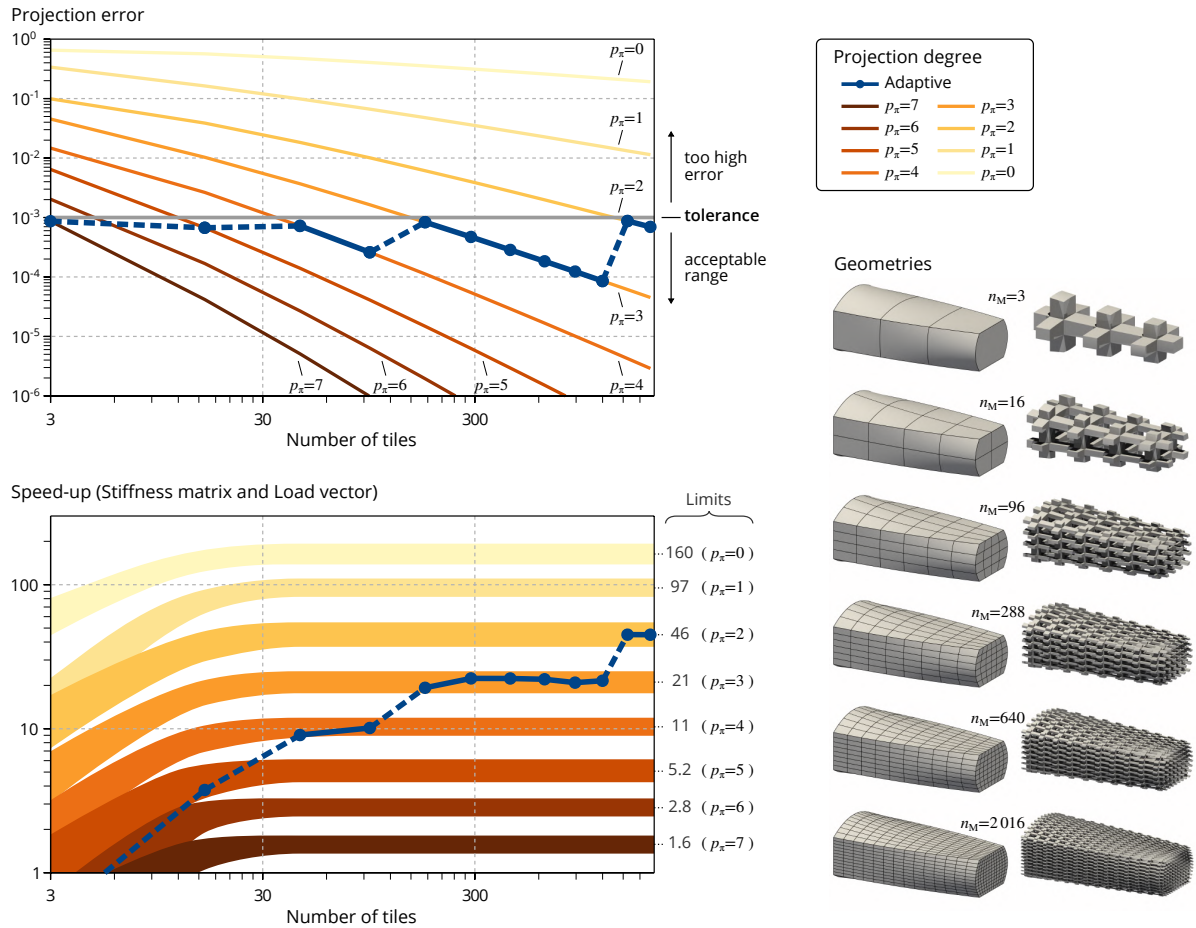

**Fig. 9** Influence of the number of tiles for a given macro-geometry: Employing the adaptive selection of the projection degree gains access to the full potential of the multiscale matrix formation. The required projection degree decreases with the increase of the number of tiles (because the projection errors due). The assembly time in comparison with the standard procedure becomes more and more appealing as reflected by the increase of the speed-up.

By comparing the cases marked by the identifiers G1, G3, and G7, one can observe the influence of the mesh density $n_T$. The degree of the solution space is the same in these three cases (quadratic). The computational time, with either the standard or the fast procedures, increases with the refinement of the solution space which is not surprising. However, the increase for the fast procedure is slower than for the standard procedure. Thus, the speed-up gets bigger while refining the mesh. The fast procedure offers the possibility to handle finer mesh sizes in an acceptable amount of time: The computation of the non-zero components of the stiffness matrix is reduced from 3 minutes with the standard procedure to 6 seconds with the presented multiscale approach.

Finally, the complexity of the shape of the reference tiles does not play a major role on the computational time. The cases marked by the identifiers G5 and G6 require the same times to form the matrices.

3.3.3 Influence of the density of the macro-geometry

Figure 9 summarizes the study of the influence of the macro-mesh density, i.e. the number of tiles. In this study, the reference tile is fixed and taken as the case G1 depicted in Figure 8. The macro-geometry is successively refined such that the number of tiles is increased by one in every direction at each step. If the projection degree is fixed then the speed-up $T_{cost}$ rapidly reaches a plateau (for $m_M \geq 30$). The computational cost increases similarly for both formation proce-dures. For instance, the speed-up is about $T_{\text{cost}} = 11 (\pm 1)$ when the projection degree is $p_\pi = 4$ when the number of tiles is $m_M \geq 30$ (see again Figure 9). However, the projection error for a given projection degree decreases with the refinement of the macro-geometry as shown in Figure 9. When adopting the strategy defined by Equation (65), the required projection degree tends to decrease while increasing the number of tiles. For instance, the required projection degree is $p_\pi = 4$ when the number of tiles is $m_M = 45$, whereas $p_\pi = 2$ when $m_M = 1573$. Consequently, the computational time for assembling the operators is kept low when following the adaptive selection of the projection degree. All the cases, from $m_M = 3$ to $m_M = 2016$, require less than 4 seconds to compute the non-zero components of the operators. Interestingly, the maximal time which is 3.5 seconds concerns the case $m_M = 1200$ for which cubic projection degree is required. The finest case, i.e. $m_M = 2016$, requires less than 3 seconds to generate the non-zero data. The speed-up is about $T_{\text{cost}} = 45$ for this last case, which attests the drastic reduction of the computational time in comparison with the standard procedure where a bit more than 2 minutes are needed to form the operators.

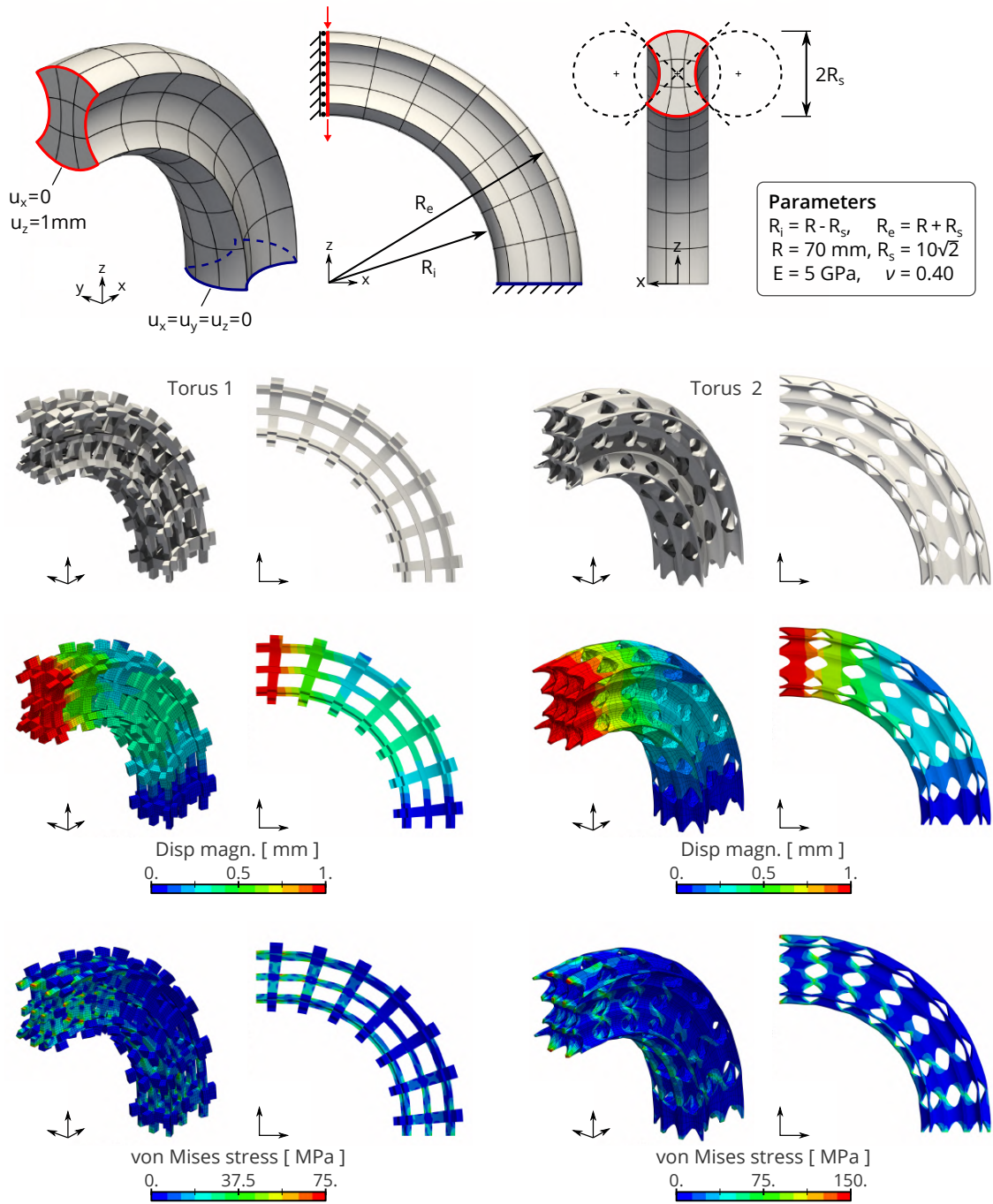

**Fig. 10** The microstructured torus-like shape subjected to an imposed displacement: definition of the problem and structural analysis results.

|  | Projection | | Numerical solutions | | | | Assembly time | | |
|---|---|---|---|---|---|---|---|---|---|
|  | $E_{\text{proj}}$ | $p_\pi$ | $E_{\text{disp}}$ | $E_{\text{stress}}$ | $E_{\text{mat}}$ | $E_{\text{vec}}$ | Gauss | Fast | $T_{\text{cost}}$ |
| Torus 1 | $7.2 \times 10^{-4}$ | 4, 4, 4 | $8.1 \times 10^{-5}$ | $9.3 \times 10^{-5}$ | $5.8 \times 10^{-5}$ | $4.3 \times 10^{-6}$ | 5.3min | 4.4s | 71 |
| Torus 2 | $7.2 \times 10^{-4}$ | 4, 4, 4 | $2.0 \times 10^{-4}$ | $1.6 \times 10^{-4}$ | $8.1 \times 10^{-5}$ | $1.5 \times 10^{-4}$ | 4.0min | 4.7s | 51 |
| Twist 1 | $1.9 \times 10^{-4}$ | 3, 3, 4 | $3.9 \times 10^{-3}$ | $3.6 \times 10^{-3}$ | $2.5 \times 10^{-5}$ | $6.2 \times 10^{-6}$ | 1.1min | 3.2s | 36 |
| Twist 2 | $1.9 \times 10^{-4}$ | 3, 3, 4 | $6.5 \times 10^{-3}$ | $3.9 \times 10^{-3}$ | $2.6 \times 10^{-5}$ | $8.5 \times 10^{-4}$ | 11.min | 6.1s | 107 |

**Table 1** Quantitative results associated to the problems depicted in figures 10 and 11.

3.3.4 Influence of the complexity of the macro-geometry

Finally, Figures 10 and 11 present additionally test cases where the macro-geometries are more challenging. The macro-geometry in Figure 10 consists in a quarter torus with a cross-section defined by the combination of three circles. This geometry is parameterized exactly with a quadratic NURBS. This macro-geometry is discretized into 63 elements. Using the tolerance on the projection error as tol $= 10^{-3}$ leads to select the projection degree as $p_\pi = 4$ (which implies $E_{\text{proj}} = 7.2 \cdot 10^{-4} \leq$ tol). Figure 10 also presents the mechanical settings of the problem. Two reference microstructures are considered here: a straight cross tile similar to the case G3 in Figure 8 but discretized with cubic splines, and the reference microstructure with identifier G8 in Figure 8. Figure 10 presents several results from the structural analyses for both microstructured-geometries. The depicted results are obtained with the operators built with the presented multiscale procedure but very similar results are obtained with the standard operators. To highlight the similarities between these numerical solutions, we provide in Table 1 the associated errors. This table also provides the computational time required to form the operators. It can be seen that the new procedure leads to the same results than the standard procedure within a much smaller amount of time.

The macro-geometry in Figure 11 consists in a helix volume (twisted extrusion of a square). Interestingly, the complexity of the geometry concerns the direction of the extrusion only. The squared cross-section can be represented with linear polynomials whereas higher degree is required to represent (approximately) the twist. Here, we employ a NURBS with degree one in the two cross-section directions and degree two in the direction of the extrusion. In order to take into account these geometrical differences during the assembly, it seems natural to select a different projection degree in the cross-section directions than in the direction of the extrusion. This can be done by initializing the projection degrees by the ones of the underlying macro-geometry. Then, the projection space is enriched (or coarsened) until the projection error is below the chosen tolerance (strategy described in Equation (65)). For the macro-geometry of the twist as given in Figure 11, it leads to a projection space with degrees $p_\pi = 3, 3, 4$. Figure 11 also describes the mechanical problem: the microstructured-geometries are subjected to their own weights. We use the same projection space for approximating the body forces (recall Equation (40)) than the one used for the stiffness matrix. Again, two reference microstructures are considered here: the reference microstructures with identifier G4 and G9 in Figure 8. The observations are the same than for the previous example of the torus-like microstructured geometry. The new multi-scale approach requires much less computational time for the formation of the finite element operators than the standard procedure, and it leads to (visually) indistinguishable numerical results as evidenced by the results in table 1. The solution errors are a bit higher than for the previous examples but are still acceptably low.

Finally, let us point out that the new procedure is successful due to the smoothness of the fields to be projected. Even for complex macro-geometries like the torus or the twist, relatively low projection degrees were found appropriate to get satisfying numerical results. Surely, this is no longer valid in case of geometrical singularities. When a macro-element is degenerated, one would preferably used more standard procedure to form the local operators there.

# 4 Extension to other quantities of interest

The presented fast assembling strategy can be extended to the computation of other quantities of interest than the finite element operators. We show briefly in this section how the developed multi-scale approach enables to perform, for instance, sensitivity analyses commonly involved in design optimization. We also discuss the construction of a matrix-free solver for the resolution.

## 4.1 Sensitivity Analysis

The volume and the compliance are often involved in design optimization problems. It is possible to perform the sensitivity analysis of these response functions with respect to the macro-geometry which follows the multiscale philosophy presented in this work.

18         T. Hirschler et al.

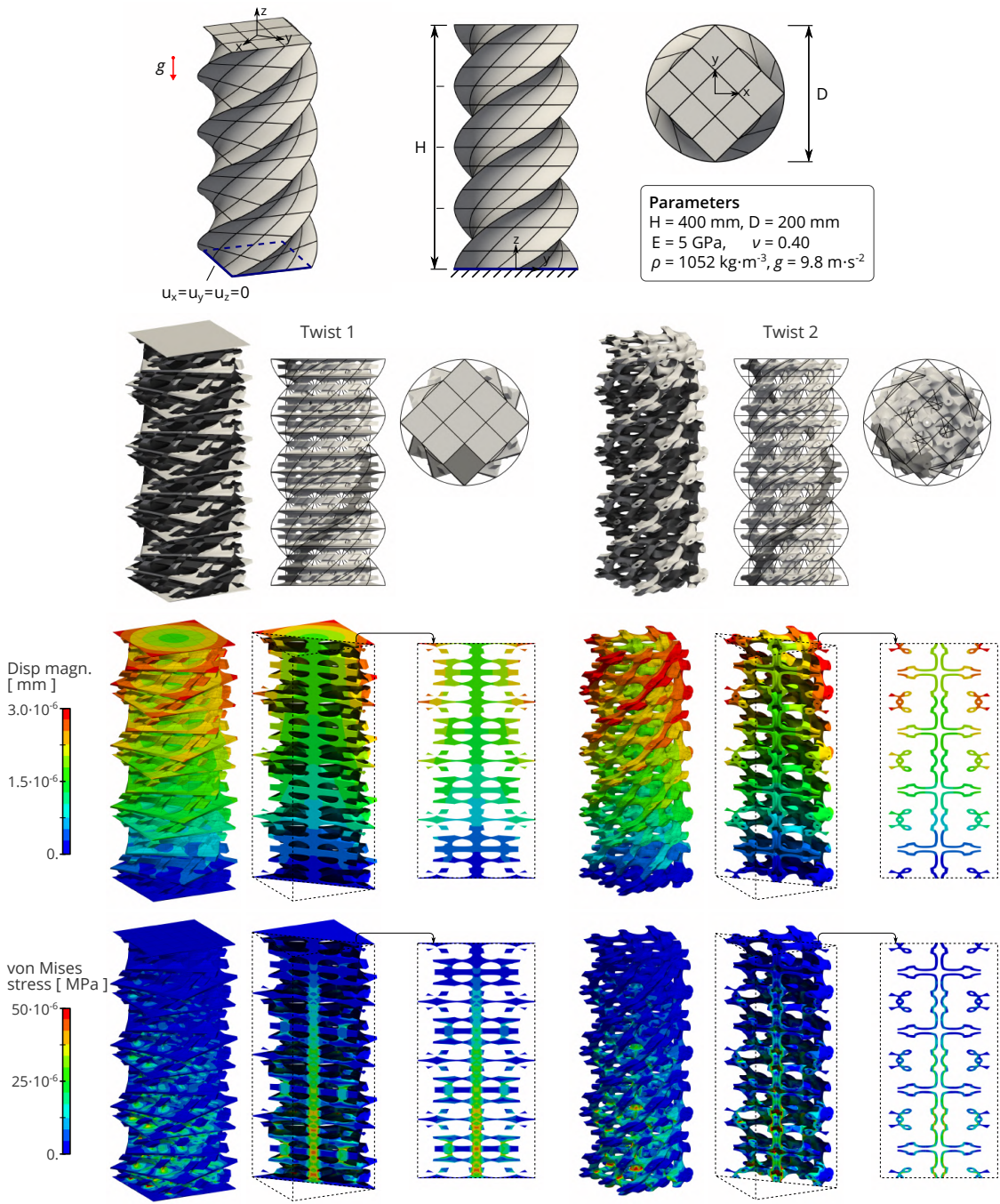

**Fig. 11** The microstructured twisted-square subjected to earth's gravity: definition of the problem and structural analysis results.

### 4.1.1 Geometric quantities

Again, the case of the volume enables to highlight the methodology. The goal is to evaluate the influence of the macro-geometry on the volume. More specifically, we seek to compute the gradient of the volume with respect to the control point coordinates associated to the macro-mapping. Then, this gradient can be employed in gradient-based optimization algorithm to solve design optimization problems.

After the split of the scale (recall Figure 4), all the geometric information associated to the macro-geometry is concatenated into the projection coefficients. Thus, the differentiation of the volume w.r.t. the macro-geometry concerns these projection coefficients principally. All the other quantities involved in the expression of the volume remain untouched (especially the lookup table) and can be concatenated into what can be seen as an adjoint field. Mathematically, the mentioned

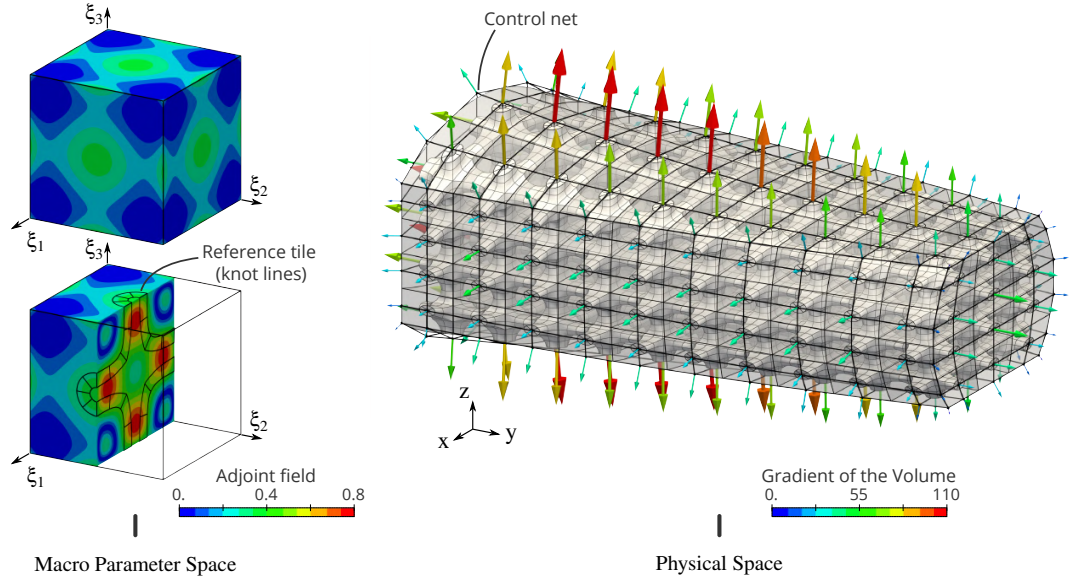

**Fig. 12** Computation of the gradient of the volume with respect to the control point of the macro-geometry. The adjoint field contains all the geometrical information regarding the reference micro-structure. The quartic polynomial degree of the projection space is considered here.

derivatives w.r.t. a control point **m** are expressed as:

$$\frac{dV^{\pi(t)}}{d\mathbf{m}} = \int_{[0,1]^{\tilde{d}}} v^{\mathrm{r}}(\xi) \frac{d}{d\mathbf{m}} |\det J_{\mathcal{M}^{(t)}}|(\xi) \, d\xi, \quad (66)$$

where $v^{\mathrm{r}}$ denotes the adjoint field, that lives in the projection space (21) and encapsulates all the information regarding the micro-structure through the Bezier-coefficients:

$$v^{\mathrm{r}}(\xi) = \sum_{C=1}^{n_{\pi}} N_C^H(\xi) \langle \mathbf{t}_{\pi}^h \rangle_C. \quad (67)$$

Here, the coefficients are basically the components of the lookup table where the inverse of the projection matrix was included (see Equation (26)). Interestingly, when taking $p_{\pi} = 0$, the adjoint field is constant and equal to the volume fraction filled by the reference microstructure.

Except this adjoint field, the differentiation of the Jacobian determinant w.r.t. the control point coordinates in Equation (66) is a standard mathematical operation found also in the classical isogeometric shape optimization framework, see for example (38) for the calculation details. Finally, the computational cost of the gradient as given by Equation (66) is almost the same than computing the gradient associated to the macro-geometry without the microstructure. The speed-up is therefore important (similar to the case of the volume itself) in comparison with the standard computation based on full quadrature rule occurring at the micro-scale. Figure 12 gives the results of the sensitivity analysis of the volume w.r.t. the control point of the macro-geometry for the same microstructured geometry than the one already studied in Section 3.1.

#### 4.1.2 Compliance

The compliance characterizes the overall stiffness of the structure and how it deforms under the prescribed loads. In design optimization, this quantity is often used as the objective function to be minimized, see for example (60). It is basically the deformation energy and can be expressed as:

$$C = \frac{1}{2} W_{\mathrm{ext}}. \quad (68)$$

Then, the derivative of the local compliance w.r.t. the control point coordinates associated to the macro-geometry can be expressed as:

$$\frac{dC^{\pi(t)}}{d\mathbf{m}} = \frac{\partial W_{\mathrm{ext}}^{\pi(t)}}{\partial \mathbf{m}} + \frac{1}{2} \frac{\partial W_{\mathrm{int}}^{\pi(t)}}{\partial \mathbf{m}}. \quad (69)$$

The remaining derivatives take the same form than the ones introduced for the case of the volume, *i.e.* Equation (66). They involve some adjoint fields and the derivatives of the macro-fields involved in the formulation of the virtual works (*i.e.* Equations (37), (40), and (41)). For instance, the differentiation of the local internal work reads as:

$$\frac{\partial W_{\mathrm{int}}^{\pi(t)}}{\partial \mathbf{m}} = \int_{[0,1]^{\tilde{d}}} \left( \sum_{ijkl=1}^{d} \mathbf{W}_{\pi,ijkl}^{h(t)}(\xi) \frac{d}{d\mathbf{m}} \bar{\mathbf{A}}_{ijkl}^{H(t)}(\xi) \right) d\xi, \quad (70)$$

where the corresponding adjoint quantities live in the projection space again:

$$\mathbf{W}_{\pi,ijkl}^{h(t)}(\xi) = \sum_{C=1}^{n_{\pi}} N_C^H(\xi) \langle \mathbf{W}_{\pi,ijkl}^{h(t)} \rangle_C, \quad \xi \in [0,1]^{\tilde{d}}. \quad (71)$$

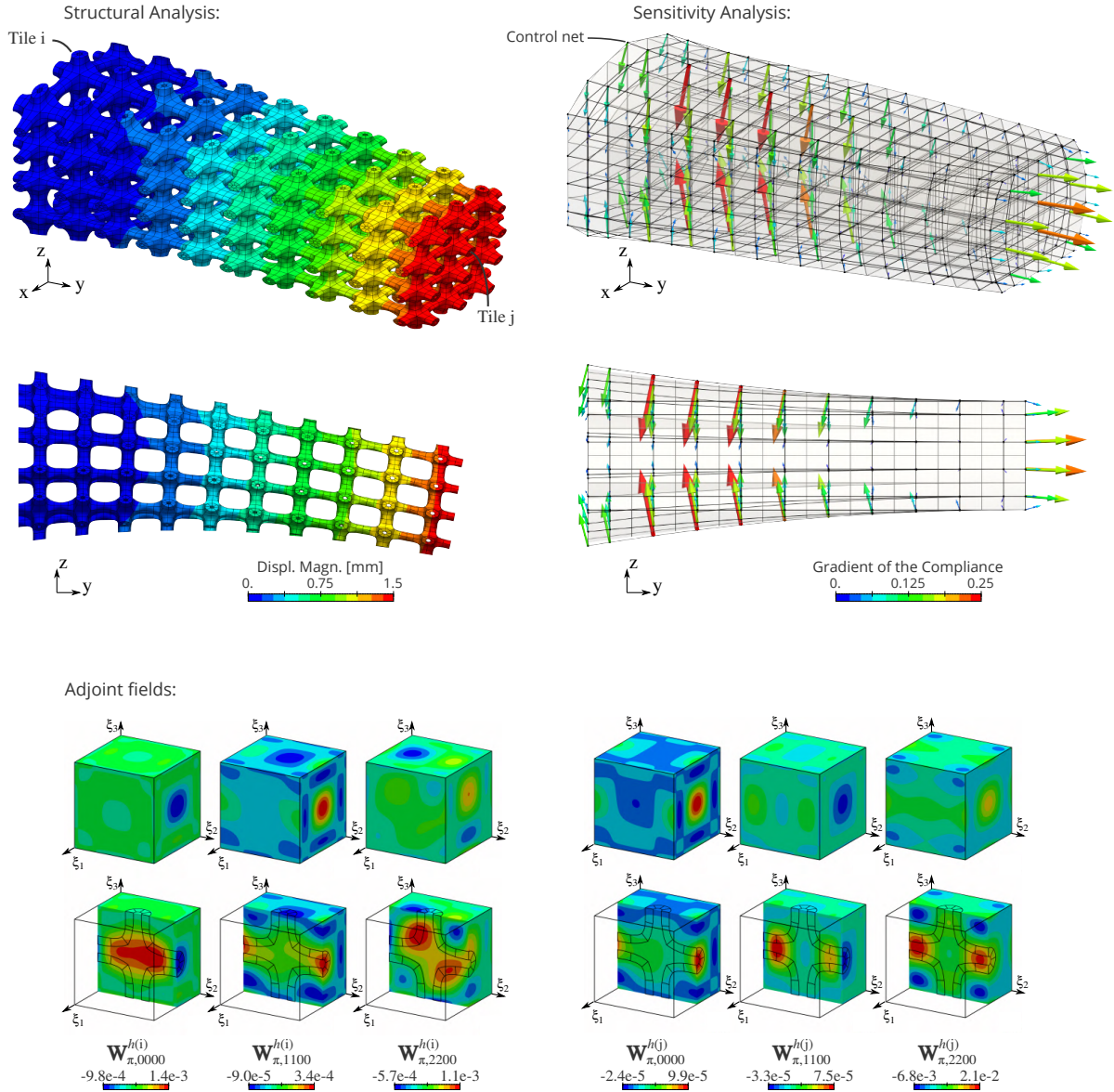

**Fig. 13** Computation of the gradient of the compliance with respect to the control points of the macro-geometry. The mechanical problem is the same than in Figure 7.

The associated Bezier coefficients are obtained by solving linear systems of the form:

$$\langle \mathbf{W}^{h(t)}_{\pi,ijkl} \rangle = \mathbf{M}^{H\,-1}_{\pi} \langle \mathbf{W}^{h(t)}_{ijkl} \rangle, \qquad (72)$$

where the right-hand sides involve the local state solution (13):

$$\langle \mathbf{W}^{h(t)}_{ijkl} \rangle_C = \int_{\Omega_T} u^{(t)}_{k,i} u^{(t)}_{l,j} N^H_C \, \mathrm{d}\xi. \qquad (73)$$

These integrals are computed via the product of the lookup table already introduced for the fast assembly of the stiffness matrix and the displacement DOF. Let us mentioned that the resolution of the linear systems as given in Equation (72), can be done off-line by integrating the inverse directly into the lookup table as done for the case of the volume (recall Equation (26) and the related discussion). Finally, the derivatives of the macro-fields w.r.t. the control points in Equation (70) are standard calculation steps in structural shape optimization, see for example (38; 60) for the details. The differentiation of the external work in Equation (69) is done similarly. Furthermore, this term might be neglected if the load support is kept identical during the design optimization (which is often the case). Employing semi-analytical differentiation schemes (10) or even automatic differentiation tools (30) might good alternatives to perform these calculation steps too.

Figure 13 shows several results regarding the sensitivity analysis for the case of the compliance. Several components of the adjoint field, for two different sub-domains, are depicted.

At the end, this gradient can be useful in an optimization process to improve the design of the structure.

### 4.2 Toward a matrix-free approach

Despite the computational time, one another important resource is the computer memory. The full matrices associated to a microstructured geometry can be large especially when the macro-geometry and/or the reference microstructure are finely discretized. The presented multiscale approach can be the starting point to design a matrix-free approach (27; 55; 61). The main idea consists in solving the linear system (30) without assembling the full matrix in a preliminary step. To that purpose, one can employ a Krylov subspace method along with a preconditioner, where only matrix-vector multiplications are performed. Interestingly, it is possible to perform the multiplication of a vector with the system matrix with only the lookup table and the projection coefficients at hand:

$$\mathbf{Ku} = \sum_{t=1}^{m_\mathrm{M}} \mathrm{LinOp}\left(\langle \bar{\mathbf{A}}^{H(t)} \rangle, \mathfrak{K}_\mathrm{r}^h, \mathbf{u}^{(t)}\right). \tag{74}$$

The linear operator LinOp computes somehow the local matrices on-the-fly. However, the local matrices are never stored but only the lookup table and the projection coefficients are. For this matrix-free approach to be valid, the preconditioning step needs to follow the same philosophy. For that purpose, the use of Domain-Decomposition Methods appears to us very attractive, see for example (12; 34). In order to push the multiscale philosophy during the resolution too, one may apply multiscale DDMs (28; 44; 68).

## 5 Conclusions

We developed a novel procedure for the formation of the finite element operators associated to microstructured geometries modeled via compositions of splines. The approach involves two principal ingredients: a polynomial approximation of the common factors in the integrals occurring at the macro-scale and lookup tables with precomputed integrals associated to the micro-scale. As a result, we exploited not only the multiscale nature of the problem, but also the repetitiveness of the micro-scale pattern during the formation of the operators. We observed major computational savings in comparison with the standard formation procedure based on Gauss quadrature. Even for low degree solution space, the assembly time is reduced from minutes (standard procedure) to seconds (new procedure) as revealed by the presented numerical examples. The speedup became overwhelming for high-order and/or refined splines. We believe the developed multiscale approach constitutes an interesting basis for the computational design of complex microstructured geometries. In order to highlight this point, we extended this fast procedure to sensitivity analyses commonly involved in design optimization. We also briefly described the construction of a matrix-free approach which consists a perspective of interest to set up an efficient resolution of the large linear systems to be solved.

*Acknowledgment.* This research was supported by the European Union Horizon 2020 research and innovation program, under grant agreement No 862025 (ADAM2), by the European Research Council through the H2020 ERC Advanced Grant 2015 No 694515 (CHANGE), and by the Swiss National Science Foundation through the project No 40B2-0_187094 (BRIDGE Discovery 2019).

## A Macro-field in linear elasticity.

We provide here several intermediary calculation steps required to identify the quantities to project during the fast-assembly via table lookup in the context of linear elasticity. The following results can also be used in standard IGA together with the methodologies developed in (49; 50; 58), for instance.

### A.1 Elastostatics with curvilinear coordinates

The virtual work (32, 33) can be formulated with respect to general curvilinear coordinates. For this purpose, one can employ the formalism commonly adopted for shell formulations, as presented in Bischoff et al. (9) for instance. The following steps consider 3-dimensional problems (*i.e.* we fix $\bar{d} = d = 3$).

By pulling back the macro-mapping $\mathcal{M}^{(t)}$, the virtual internal work is transformed as:

$$W_\mathrm{int}^{(t)}(\boldsymbol{u}, \boldsymbol{v}) = -\int_{\Omega^{(t)}} \boldsymbol{\sigma}(\boldsymbol{u}) : \boldsymbol{\varepsilon}(\boldsymbol{v}) \, \mathrm{d}\Omega, \tag{75}$$

$$= -\int_{\Omega_\mathrm{T}} \left(\boldsymbol{\sigma}^{(t)}(\boldsymbol{u}) : \boldsymbol{\varepsilon}^{(t)}(\boldsymbol{v})\right) |\det J_{\mathcal{M}^{(t)}}| \, \mathrm{d}\xi. \tag{76}$$

where the stress and strain tensors are expressed with respect to the macro-curvilinear coordinates $\xi_i$ (recall the notations in Figure 2). The strain tensors are represented by their covariant components $\hat{\varepsilon}_{ij}$ and the stress tensors are represented by their contravariant components $\check{\sigma}_{ij}$:

$$\boldsymbol{\varepsilon}^{(t)} = \hat{\varepsilon}_{ij}^{(t)} \, \check{\mathbf{g}}_i^{(t)} \otimes \check{\mathbf{g}}_j^{(t)}, \qquad \boldsymbol{\sigma}^{(t)} = \check{\sigma}_{ij}^{(t)} \, \hat{\mathbf{g}}_i^{(t)} \otimes \hat{\mathbf{g}}_j^{(t)}, \tag{77}$$

where the Einstein's summation convention applies, and $\otimes$ defines the tensor product of two vectors. With this choice, the scalar product of the strain and stress second-order tensors (operator : in the internal work) reads as:

$$\boldsymbol{\sigma}^{(t)} : \boldsymbol{\varepsilon}^{(t)} = \sum_{i=1}^{3} \sum_{j=1}^{3} \check{\sigma}_{ij}^{(t)} \, \hat{\varepsilon}_{ij}^{(t)}. \tag{78}$$

The covariant basis vectors $\hat{\mathbf{g}}_i$ are obtained from partial derivatives of the global mappings with respect to the macro-parameters, which are abbreviated via subscripted comma for short notation:

$$\hat{\mathbf{g}}_i^{(t)} = \mathcal{M}_{,\xi_i}^{(t)}. \tag{79}$$

Their orthogonal counterpart, *i.e.* the contravariant basis vectors $\check{\mathbf{g}}_i$, can be obtained through:

$$\check{\mathbf{g}}_i^{(t)} = \check{g}_{i1}^{(t)} \hat{\mathbf{g}}_1^{(t)} + \check{g}_{i2}^{(t)} \hat{\mathbf{g}}_2^{(t)} + \check{g}_{i3}^{(t)} \hat{\mathbf{g}}_3^{(t)} \tag{80}$$

where the contravariant metrics $\check{g}_{ij}$ result from inverting 3-by-3 matrices formed by the covariant metrics $\hat{g}_{ij}$:

$$\left[\check{g}_{ij}^{(t)}\right] = \left[\hat{g}_{ij}^{(t)}\right]^{-1}, \quad \text{where} \quad \hat{g}_{ij}^{(t)} = \hat{\mathbf{g}}_i^{(t)} \cdot \hat{\mathbf{g}}_j^{(t)}. \tag{81}$$

With these notations in hand, the covariant strains components in Equation (78) are given by:

$$\hat{\varepsilon}_{ij}^{(t)}(\boldsymbol{u}) = \frac{1}{2}\left(\boldsymbol{u}_{,\xi_i} \cdot \hat{\mathbf{g}}_j^{(t)} + \boldsymbol{u}_{,\xi_j} \cdot \hat{\mathbf{g}}_i^{(t)}\right). \tag{82}$$

In the case of isotropic materials with an elastic behavior, the contravariant stresses components also involved in Equation (78) are given by:

$$\check{\sigma}_{ij}^{(t)}(\boldsymbol{u}) = \sum_{k=1}^{3}\sum_{l=1}^{3} \check{C}_{ijkl}^{(t)} \hat{\varepsilon}_{kl}^{(t)}(\boldsymbol{u}), \tag{83}$$

where the contravariant components of the material tensor are:

$$\check{C}_{ijkl}^{(t)} = \lambda\, \check{g}_{ij}^{(t)} \check{g}_{kl}^{(t)} + \mu\left(\check{g}_{ik}^{(t)}\check{g}_{jl}^{(t)} + \check{g}_{il}^{(t)}\check{g}_{jk}^{(t)}\right), \tag{84a}$$

$$\check{C}_{ijkl}^{(t)} = \check{C}_{jikl}^{(t)} = \check{C}_{ijlk}^{(t)} = \check{C}_{klij}^{(t)}, \tag{84b}$$

which involve two material parameters, *i.e.* the so-called Lamé constants $\lambda$ and $\mu$. These coefficients are usually prescribed using the Young's modulus $E$ and the Poisson's ratio $\nu$ of the material.

The external virtual work also needs to be manipulated in order to apply the fast assembly strategy. Similarly to the internal work, the integrals defining the external work need to be expressed with respect to the macro-curvilinear coordinates. Performing the change of variables by pulling back the macro-mappings is quite straightforward for the case of the body forces. On the contrary, dealing with the tractions may require more attention. Two loading scenarios should be distinguished as defined in Figure 5: the case where load support overlaps the macro-boundaries, and the case where the load support lies inside the macro-elements. This last case will not be considered in the scope of this work. Consequently, after the change of variables, the considered local external works are of the form:

$$W_{\text{ext}}^{(t)}(\boldsymbol{v}) = \int_{\Omega_T} \left(\boldsymbol{v} \cdot \bar{\mathbf{b}}^{(t)}\right) d\xi + \int_{\Gamma_T} \left(\boldsymbol{v} \cdot \bar{\mathbf{t}}^{(t)}\right) d\Gamma_T. \tag{85}$$

where

$$\bar{\mathbf{b}}^{(t)}(\xi) = \left(\mathbf{b} \circ \mathcal{M}^{(t)}\right)(\xi) \, |\det J_{\mathcal{M}^{(t)}}|(\xi), \; \xi \in \Omega_T, \tag{86}$$

$$\bar{\mathbf{t}}^{(t)}(\xi) = \left(\mathbf{t} \circ \mathcal{M}^{(t)}\right)(\xi) \, |\det I_{\mathcal{M}^{(t)}}|(\xi), \; \xi \in \Gamma_T \cap \partial[0,1]^3, \tag{87}$$

and

$$\det I_{\mathcal{M}^{(t)}} = \left(\mathbf{v}_1^{(t)} \times \mathbf{v}_2^{(t)}\right) \cdot \mathbf{v}_3^{(t)},$$
$$\mathbf{v}_i^{(t)} = \begin{cases} \mathbf{n}_i^{(t)} = \check{\mathbf{g}}_i^{(t)}/|\check{\mathbf{g}}_i^{(t)}| & \text{if } \xi_i \in \{0,1\}, \\ \hat{\mathbf{g}}_i^{(t)} & \text{else.} \end{cases} \tag{88}$$

Further information regarding continuum formulation with curvilinear coordinates can be found, for instance in Bischoff et al. (9).

## A.2 Expression of the macro-fields in the internal work

### A.2.1 Matrix formulation

In order to reveal the macro-quantities to be projected, one can use the following matrix formulation. The following steps consider 3-dimensional problems (*i.e.* we fix $\bar{d} = d = 3$). Due to the symmetries in the strains, the stresses, and the material tensors, a compact matrix representation of theses quantities is possible. By employing Voigt notation, the six distinct strains components are concatenated into a vector given by:

$$\hat{\boldsymbol{\varepsilon}}^{(t)}(\boldsymbol{u}) = \hat{\mathbf{G}}^{(t)\top} \boldsymbol{u}_{,\xi}, \tag{89}$$

with $\hat{\boldsymbol{\varepsilon}}^{(t)} \in \mathbb{R}^6$, $\hat{\mathbf{G}}^{(t)\top} \in \mathbb{R}^{6\times 9}$, and $\boldsymbol{u}_{,\xi} \in \mathbb{R}^9$. According to Equation (82), their expressions are:

$$\hat{\boldsymbol{\varepsilon}}^{(t)} = \begin{pmatrix} \hat{\varepsilon}_{11} \\ \hat{\varepsilon}_{22} \\ \hat{\varepsilon}_{33} \\ 2\hat{\varepsilon}_{23} \\ 2\hat{\varepsilon}_{13} \\ 2\hat{\varepsilon}_{12} \end{pmatrix}^{(t)}, \quad \hat{\mathbf{G}}^{(t)\top} = \begin{bmatrix} \hat{\mathbf{g}}_1^\top & & \\ & \hat{\mathbf{g}}_2^\top & \\ & & \hat{\mathbf{g}}_3^\top \\ & \hat{\mathbf{g}}_3^\top & \hat{\mathbf{g}}_2^\top \\ \hat{\mathbf{g}}_3^\top & & \hat{\mathbf{g}}_1^\top \\ \hat{\mathbf{g}}_2^\top & \hat{\mathbf{g}}_1^\top & \end{bmatrix}^{(t)}, \quad \boldsymbol{u}_{,\xi} = \begin{pmatrix} \boldsymbol{u}_{,\xi_1} \\ \boldsymbol{u}_{,\xi_2} \\ \boldsymbol{u}_{,\xi_3} \end{pmatrix}. \tag{90}$$

Key sub-blocks denoted $\hat{\mathbf{G}}_i^{(t)} \in \mathbb{R}^{3\times 6}$ can be identified in the strains-displacement relations:

$$\hat{\mathbf{G}}_1^{(t)} = \begin{bmatrix} \hat{\mathbf{g}}_1 & \mathbf{0} & \mathbf{0} & \mathbf{0} & \hat{\mathbf{g}}_3 & \hat{\mathbf{g}}_2 \end{bmatrix}^{(t)}, \tag{91a}$$

$$\hat{\mathbf{G}}_2^{(t)} = \begin{bmatrix} \mathbf{0} & \hat{\mathbf{g}}_2 & \mathbf{0} & \hat{\mathbf{g}}_3 & \mathbf{0} & \hat{\mathbf{g}}_1 \end{bmatrix}^{(t)}, \tag{91b}$$

$$\hat{\mathbf{G}}_3^{(t)} = \begin{bmatrix} \mathbf{0} & \mathbf{0} & \hat{\mathbf{g}}_3 & \hat{\mathbf{g}}_2 & \hat{\mathbf{g}}_1 & \mathbf{0} \end{bmatrix}^{(t)}, \tag{91c}$$

such that:

$$\hat{\mathbf{G}}^{(t)\top} = \begin{bmatrix} \hat{\mathbf{G}}_1^{(t)\top} & \hat{\mathbf{G}}_2^{(t)\top} & \hat{\mathbf{G}}_3^{(t)\top} \end{bmatrix}. \tag{92}$$

Then, the strains-stresses relation (the constitutive law, already introduced in Equation (83)) reads, in the matrix form, as:

$$\begin{pmatrix} \check{\sigma}_{11} \\ \check{\sigma}_{22} \\ \check{\sigma}_{33} \\ \check{\sigma}_{23} \\ \check{\sigma}_{13} \\ \check{\sigma}_{12} \end{pmatrix}^{(t)} = \begin{bmatrix} \check{C}_{1111} & \check{C}_{1122} & \check{C}_{1133} & \check{C}_{1123} & \check{C}_{1113} & \check{C}_{1112} \\ * & \check{C}_{2222} & \check{C}_{2233} & \check{C}_{2223} & \check{C}_{2213} & \check{C}_{2212} \\ * & * & \check{C}_{3333} & \check{C}_{3323} & \check{C}_{3313} & \check{C}_{3312} \\ * & * & * & \check{C}_{2323} & \check{C}_{2313} & \check{C}_{2312} \\ * & * & * & * & \check{C}_{1313} & \check{C}_{1312} \\ * & * & * & * & * & \check{C}_{1212} \end{bmatrix}^{(t)} \begin{pmatrix} \hat{\varepsilon}_{11} \\ \hat{\varepsilon}_{22} \\ \hat{\varepsilon}_{33} \\ 2\hat{\varepsilon}_{23} \\ 2\hat{\varepsilon}_{13} \\ 2\hat{\varepsilon}_{12} \end{pmatrix}^{(t)}, \tag{93}$$

where the given matrix, later denoted $\check{\mathbf{C}}^{(t)}$, is symmetric and contains the 21 distinct components of the material tensor. By substituting Equation (89) in Equation (93), the stresses can be expressed with respect to the displacement field via:

$$\check{\boldsymbol{\sigma}}^{(t)}(\boldsymbol{u}) = \check{\mathbf{C}}^{(t)} \hat{\boldsymbol{\varepsilon}}^{(t)}(\boldsymbol{u}) = \check{\mathbf{C}}^{(t)} \hat{\mathbf{G}}^{(t)\top} \boldsymbol{u}_{,\boldsymbol{\xi}}. \tag{94}$$

Finally, this leads us to the following expression of the internal virtual work:

$$W_{\text{int}}^{(t)}(\boldsymbol{u}, \boldsymbol{v}) = -\int_{\Omega_T} \left( \check{\boldsymbol{\sigma}}^{(t)}(\boldsymbol{u}) \cdot \hat{\boldsymbol{\varepsilon}}^{(t)}(\boldsymbol{v}) \right) |\det J_{\mathcal{M}^{(t)}}| \, \mathrm{d}\boldsymbol{\xi}, \tag{95}$$

$$= -\int_{\Omega_T} \left( \boldsymbol{u}_{,\boldsymbol{\xi}}^\top \hat{\mathbf{G}}^{(t)} \check{\mathbf{C}}^{(t)} \hat{\mathbf{G}}^{(t)\top} \boldsymbol{v}_{,\boldsymbol{\xi}} \right) |\det J_{\mathcal{M}^{(t)}}| \, \mathrm{d}\boldsymbol{\xi}, \tag{96}$$

$$= -\sum_{i=1}^{3} \sum_{j=1}^{3} \int_{\Omega_T} \left( \boldsymbol{u}_{,\xi_i}^\top \hat{\mathbf{G}}_i^{(t)} \check{\mathbf{C}}^{(t)} \hat{\mathbf{G}}_j^{(t)\top} \boldsymbol{v}_{,\xi_j} \right) |\det J_{\mathcal{M}^{(t)}}| \, \mathrm{d}\boldsymbol{\xi}. \tag{97}$$

This last equation enables to identify the macro-fields to project as provided in Equation (36).

### A.2.2 A simple case for validation

For illustrative purpose (and *code* verification), the macro-field involved in the internal work can be explicitly expressed for simple geometries. Let us for instance consider the macro-geometry as a cube of length $L$; *i.e.* the global mapping reads:

$$\mathcal{M}(\xi_1, \xi_2, \xi_3) = (\xi_1 L)\,\mathbf{x} + (\xi_2 L)\,\mathbf{y} + (\xi_3 L)\,\mathbf{z}, \tag{98}$$

with $\xi_1, \xi_2, \xi_3 \in [0, 1]$. In this simple case, the macro-fields are constants, and are given by:

$$\bar{\mathbf{A}}_{11}^H = \begin{bmatrix} \alpha_1 & & \\ & \alpha_2 & \\ & & \alpha_2 \end{bmatrix}, \quad \bar{\mathbf{A}}_{22}^H = \begin{bmatrix} \alpha_2 & & \\ & \alpha_1 & \\ & & \alpha_2 \end{bmatrix}, \quad \bar{\mathbf{A}}_{33}^H = \begin{bmatrix} \alpha_2 & & \\ & \alpha_2 & \\ & & \alpha_1 \end{bmatrix},$$
$$\bar{\mathbf{A}}_{12}^H = \begin{bmatrix} & \alpha_3 & \\ \alpha_2 & & \\ & & \end{bmatrix}, \quad \bar{\mathbf{A}}_{13}^H = \begin{bmatrix} & & \alpha_3 \\ & & \\ \alpha_2 & & \end{bmatrix}, \quad \bar{\mathbf{A}}_{23}^H = \begin{bmatrix} & & \\ & & \alpha_3 \\ & \alpha_2 & \end{bmatrix}. \tag{99}$$

where the three distinct coefficients are:

$$\alpha_1 = (\lambda + 2\mu)L, \quad \alpha_2 = \mu L, \quad \alpha_3 = \lambda L. \tag{100}$$